\documentclass[11pt]{amsart}
\usepackage[T1]{fontenc}
\usepackage{amssymb,amsmath,amsfonts,mathtools,enumitem,bm}
\usepackage[margin=1in, headsep=18pt]{geometry}
\usepackage{tikz-cd}
\usepackage{graphicx}
\usepackage{xcolor}
\usepackage[colorlinks=true,linkcolor=blue,citecolor=blue]{hyperref}
\newcommand\btitle[1]{\textit{``#1''}}
\makeatletter
\newcommand\eprint[1]{\@eprint#1\@nil}            
\def\@eprint#1:#2\@nil{\href{https://arxiv.org/abs/#2}{arXiv:#2}}
\makeatother

\numberwithin{equation}{section}
\newtheorem{theorem}[equation]{Theorem}
\newtheorem{proposition}[equation]{Proposition}
\newtheorem{lemma}[equation]{Lemma}
\newtheorem{corollary}[equation]{Corollary}
\newtheorem{conjecture}[equation]{Conjecture}
\newtheorem{introthm}{Theorem}

\theoremstyle{definition}
\newtheorem{definition}[equation]{Definition}
\newtheorem{remark}[equation]{Remark}
\newtheorem{example}[equation]{Example}
\newtheorem{notation}[equation]{Notation}

\usepackage{ulem} 
\usepackage{color}

\newcommand\cbl[1]{{\color{blue} #1}}

\newcommand{\gl}{\mathfrak{gl}}
\newcommand{\fg}{\mathfrak{g}}
\newcommand{\fh}{\mathfrak{h}}
\newcommand{\fn}{\mathfrak{n}}
\newcommand{\fb}{\mathfrak{b}}
\newcommand{\U}{\mathrm{U}}
\newcommand{\bW}{\mathbb{W}}          
\newcommand{\bA}{\mathbb{A}}          
\newcommand{\bU}{\mathbb{U}}
\newcommand{\GPD}{\mathrm{GPD}}         
\newcommand{\GH}{\mathrm{GH}}         
\newcommand{\eps}{\varepsilon}
\newcommand{\ch}{\mathrm{ch}}

\newcommand{\Hom}{\mathrm{Hom}}
\newcommand{\Ext}{\mathrm{Ext}}
\newcommand{\Ann}{\mathrm{Ann}}

\newcommand{\Gr}{\mathrm{Gr}}
\newcommand{\Rep}{\mathsf{Rep}}
\newcommand{\id}{\mathrm{id}}

\newcommand{\bx}{\bar{x}}
\newcommand{\by}{\bar{y}}
\newcommand{\Lam}{\Lambda}

\newcommand{\cJ}{\mathcal{J}}

\newcommand{\sS}{\mathfrak{S}}
\newcommand{\qbinom}[2]{\genfrac{[}{]}{0pt}{}{#1}{#2}_{q}}
\DeclareMathOperator{\Span}{span}

\begin{document}
	
	\title[Two-row Garsia--Haiman modules]{Generalized parabolic Demazure modules and \\ the $N!/2$-conjecture for two-row diagrams}
	
	\author{Anton Khoroshkin}
	\address{Anton Khoroshkin: Department of Mathematics, University of Haifa, Mount Carmel, 3498838, Haifa, Israel}
	\email{khoroshkin@gmail.com}
	
	\author{Ievgen Makedonskyi}
	\address{Ievgen Makedonskyi: Yanqi Lake Beijing Institute of Mathematical Sciences And Applications (BIMSA), No.~544, Hefangkou Village, Huaibei Town, Huairou District, Beijing 101408}
	\email{mak@bimsa.cn}
	
	
	\begin{abstract}
		We study a family of generalized parabolic Demazure modules, i.e.\ of $\gl[x]$-submodules of tensor products of integrable representations of the affine Lie algebra generated by tensor products of extremal vectors, and describe these modules by generators and relations. Realizing the Garsia--Haiman modules of two-row diagrams as fusion products of two local Weyl modules, we show that their $y$-graded components admit filtrations whose associated graded modules are isomorphic to tensor products of these generalized parabolic Demazure modules with the cohomology rings of Grassmannians. As a consequence, the Garsia--Haiman modules of two-row diagrams of a fixed size generate a distributive lattice, and the $N!/2$-conjecture holds for two-row diagrams.
	\end{abstract}

	\maketitle
	
	\setcounter{tocdepth}{1}
	\tableofcontents
	
	\section{Introduction}\label{sec:intro}
	
	This is the second paper in which we use the iterated fusion product of Feigin and Loktev~\cite{FeLo} to study Garsia--Haiman modules. In~\cite{KM} we showed that, after a Schur--Weyl duality for the current algebra $\gl_V[x,y]$, the Garsia--Haiman module of a Young diagram becomes an iterated fusion product of the local Weyl modules attached to its rows; this gave a representation-theoretic proof of Butler's conjecture and lower bounds in the $\frac{N!}{k}$-conjecture. In the present paper we work out completely the simplest nontrivial case, the diagrams with two rows. We describe the Garsia--Haiman modules of two-row diagrams by generators and relations and deduce that they generate a distributive lattice and that the $\frac{N!}{2}$-conjecture holds for them.
	
	Recall that the Garsia--Haiman module $\GH^{\mathfrak S}_\lambda$ of a Young diagram $\lambda\vdash N$ with cells $(a_c,b_c)$ is the span of all partial derivatives of the polynomial $\Delta_\lambda=\det\big(x_i^{a_c}y_i^{b_c}\big)_{i,c=1}^N$~\cite{GH}. By Haiman's theorem~\cite{Haiman} $\dim\GH^{\mathfrak S}_\lambda=N!$, and the bigraded Frobenius characteristic of $\GH^{\mathfrak S}_\lambda$ is the modified Macdonald polynomial $\widetilde H_\lambda(X;q,t)$~\cite{Hai99,Mac88}. Garsia and Haiman~\cite{GH96} conjectured that if $\nu_1,\nu_2\vdash N$ are obtained from a diagram $\lambda\vdash N+1$ by removing two different corners, then
	\[
	\dim\big(\GH^{\mathfrak S}_{\nu_1}\cap\GH^{\mathfrak S}_{\nu_2}\big)=\frac{N!}{2}
	\]
	(the $\frac{N!}2$-conjecture; see~\cite{BG} for the $\frac{N!}k$-conjecture and~\cite{But,GKX,KM} for Butler's conjecture on the character of the intersection). The $\frac{N!}2$-conjecture is known for hooks~\cite{BH,Arm}; in general only the lower bound is known~\cite{KM}.
	
	The Schur--Weyl duality of~\cite[\S3]{KM} turns $\GH^{\mathfrak S}_\lambda$ into the cyclic $\gl_V[x,y]$-module $\GH^{\gl}_\lambda=S^N(V[x,y])/\cJ_\lambda$, and by~\cite[Theorem~A]{KM} and Haiman's theorem (see Theorem~\ref{thm:GH-fusion})
	\[
	\GH^{\gl}_\lambda\ \simeq\ W_{\lambda_0\eps_1}(x)\star_yW_{\lambda_1\eps_1}(x)\star_y\cdots\star_yW_{\lambda_l\eps_1}(x).
	\]
	For a two-row diagram $(N-j,j)$ this is the $y$-fusion product $\mathcal H_j=W_{(N-j)\eps_1}\star_yW_{j\eps_1}$ of two local Weyl modules of $\gl_V[x]$. It splits into the $y$-graded components $\mathcal H_j^k$, $0\le k\le j$, which are graded $\gl_V[x]$-modules, and our results describe these components.
	
	\subsection{Generalized parabolic Demazure modules}
	
	Our main technical tool is a description by generators and relations of a family of generalized parabolic Demazure modules. Given dominant weights $\mu^1,\dots,\mu^p$ of a simply laced Lie algebra $\fg$, the module $\GPD(\mu^1,\dots,\mu^p)$ is the $\fg[x]$-submodule of $W_{\mu^1}\otimes W_{\mu^1+\mu^2}\otimes\dots\otimes W_{\mu^1+\dots+\mu^p}$ generated by the tensor product of the cyclic vectors (Definition~\ref{def:GPD}). For $p=1$ this is the local Weyl module $W_{\mu^1}$. In general these are the generalized Demazure modules of Naoi~\cite{Naoi} attached to anti-dominant extremal weights of level one; they are stable under the parabolic subalgebra $\fg[x]$, whence the name. By Dumanski~\cite{D} they are the fusion products $D(p,pw_0(\mu^1))\star D(p-1,(p-1)w_0(\mu^2))\star\dots\star D(1,w_0(\mu^p))$ of Demazure modules of different levels, where $w_0$ is the longest element of the finite Weyl group. In particular, their characters as $\fg$-modules are known.
	
	However, to the best of our knowledge, neither the defining relations nor the graded characters of these modules are known in general. Naoi's realization of generalized Demazure modules by successive applications of Demazure functors~\cite[\S2]{Naoi} expresses their graded characters through iterated Demazure operators, $D_{w_1}\big(e^{\Lambda_1}\cdot D_{w_2}\big(e^{\Lambda_2}\cdots D_{w_p}(e^{\Lambda_p})\cdots\big)\big)$~\cite{LLM}, which is far from being explicit. In Section~\ref{sec:GPD} we find the defining relations in the case $p=2$. Let $\fg=\gl_m$, $f_i=E_{i+1,i}$, and let $w$ denote the cyclic vector of a local Weyl module.
	
	\begin{introthm}[Theorem~\ref{thm:SES}, Corollaries~\ref{cor:monomial} and~\ref{cor:two-rows-presentation}]\label{thm:introA}
		\begin{enumerate}
			\item For dominant weights $\mu^1,\dots,\mu^p$, $p\ge2$, and $1\le k<m$ there is a short exact sequence of graded $\fg[x]$-modules
			\[
			0\to\GPD(\mu^1,\dots,\mu^p+\omega_{k-1}+\omega_{k+1})\to\GPD(\mu^1,\dots,\mu^p+2\omega_k)\to\GPD(\mu^1,\dots,\mu^{p-1}+\omega_k,\mu^p)\to0 .
			\]
			\item For dominant weights $\mu=\sum_im_i\omega_i$ and $\nu$, the module $\GPD(\mu,\nu)$ is the quotient of the local Weyl module $W_{2\mu+\nu}$ by the monomial relations $f_ix^{\langle\mu+\nu,\alpha_i^\vee\rangle}w=0$, $m_i>0$.
			\item In particular, for $0\le k\le j\le N/2$ the module $\GPD_{j,k}:=\GPD\big((j-k)\omega_1,(N-2j)\omega_1+k\omega_2\big)$ is the quotient of $W_{(N-k)\eps_1+k\eps_2}$ by the single relation $E_{21}x^{N-j-k}w=0$.
		\end{enumerate}
	\end{introthm}
	
	We also deduce from it an Ext-vanishing property of the modules $\GPD_{j,k}$ (Proposition~\ref{prop:Ext}), which restricts their graded characters to the span of the plethystic Schur functions 
	$$\ch_q(\GPD_{j,k})\in \mathbb{Z}\Bigl\langle s_\kappa\big[\frac X{1-q}\big] \colon \kappa\succeq(N-j,j)\Bigr\rangle.$$
	
	\subsection{The \texorpdfstring{$y$}{y}-graded components of the two-row modules}
	
	Fix $0\le k\le j\le N/2$. Let $\bW_{(N-k)\eps_1+k\eps_2}$ be the global Weyl module of $\gl_V[x]$ with the cyclic vector $w$, and let $v_k\in\mathcal H_j^k$ be the class of $(E_{21}y)^k$ applied to the cyclic vector of $\mathcal H_j$. For $c=1,2$ let $\psi e_r(t^c_\bullet)\in\U(\fh[x])$ be the elementary symmetric functions in the Cartan currents $E_{cc}x^t$ regarded as power sums:
	\[
	\sum_{r\ge0}\psi e_r(t^c_\bullet)\,z^r=\exp\Big(\sum_{t\ge1}(-1)^{t-1}\frac{E_{cc}x^t\,z^t}{t}\Big).
	\]
	
	\begin{introthm}[Theorem~\ref{thm:presentation}, Corollaries~\ref{cor:A} and~\ref{thm:B}]\label{thm:introB}
		\begin{enumerate}
			\item The assignment $w\mapsto v_k$ defines a surjection $\bW_{(N-k)\eps_1+k\eps_2}\twoheadrightarrow\mathcal H_j^k$, whose kernel is generated by
			\begin{itemize}
				\item[(R1)] $E_{cc}x^l w$ $(c\ge3)$, $(E_{11}+E_{22})x^l w$, $\psi e_r(t^2_\bullet)w$ $(r>k)$, $\psi e_s(t^1_\bullet)w$ $(s>j-k)$, where $l\ge 1$;
				\item[(R2)] $\psi e_{j-k}(t^1_\bullet)\,E_{21}x^{N-2j}\,a\,w$, where $a\in\U(\fh[x])$.
			\end{itemize}
			\item The algebra $\U(\fh[x])/\Ann(v_k)$ of the Cartan currents acting on $v_k$ is isomorphic to the cohomology ring $H^*(\Gr(k,j);\mathbb C)$.
			\item $\mathcal H_j^k$ has a filtration by $\gl_V[x]$-submodules with $\binom jk$ subquotients, each isomorphic, up to a shift of the grading, to $\GPD_{j,k}$. In particular, $\ch_q\mathcal H_j^k=\qbinom jk\ch_q\GPD_{j,k}$.
		\end{enumerate}
	\end{introthm}
	
	The relations (R1) involve only the Cartan currents: they say that the Cartan currents act on $w$ as on the weight space $W_{j\eps_1}[(j-k)\eps_1+k\eps_2]$ of the local Weyl module, which is a copy of $H^*(\Gr(k,j))$. Modulo the Cartan currents of positive degree, the relation (R2) becomes the relation $E_{21}x^{N-j-k}w=0$ of $\GPD_{j,k}$ (Theorem~\ref{thm:introA}(3)); thus (R2) is a deformation of this monomial relation by the Cartan currents. We verify the relations (R1) and (R2) directly in the fusion product, using only the relations of the local Weyl modules. We show that there are no other relations from a dimension count based on Theorem~\ref{thm:introA}.
	
	\begin{remark}\label{rem:Macdonald-intro}
		By~\cite{Hai99,Haiman} the bigraded character of $\mathcal H_j$ is the modified Macdonald polynomial $\widetilde H_{(N-j,j)}(X;q,t)$, with $q$ and $t$ recording the degrees in $x$ and $y$. Hence Theorem~\ref{thm:introB}(3) gives
		\[
		\widetilde H_{(N-j,j)}(X;q,t)\ =\ \sum_{k=0}^{j}t^k\qbinom jk\,\ch_q\GPD_{j,k}(X),
		\]
		so the graded characters of the modules $\GPD_{j,k}$ are determined by the Macdonald polynomials of two-row diagrams.
	\end{remark}
	
	\subsection{Distributivity and the \texorpdfstring{$\frac{N!}{2}$}{N!/2}-conjecture}
	
	For fixed $N$ all the modules $\mathcal H_j$ are quotients of $S^N(V[x,y])$, and Theorem~\ref{thm:introB} allows us to compare their kernels. For fixed $k$, the Cartan relations (R1) become weaker as $j$ grows, while the relation (R2) for $j$ follows from (R2) for $j+1$ together with (R1) for $j$. Hence every kernel is the sum of a member of a decreasing chain and a member of an increasing chain of submodules, and Birkhoff's lemma on two chains gives the first main result.
	
	\begin{introthm}[Theorems~\ref{thm:lattice} and~\ref{thm:lattice-gl}]\label{thm:introC}
		For every $N$ the Garsia--Haiman modules $\GH^{\mathfrak S}_{(N-j,j)}$, $0\le j\le\lfloor N/2\rfloor$, generate a distributive lattice of subspaces of $\mathbb C[x_1,y_1,\dots,x_N,y_N]$. Equivalently, the Garsia--Haiman ideals $J_{(N-j,j)}$ generate a distributive lattice of ideals, and the $\gl_V[x,y]$-submodules $\cJ_{(N-j,j)}\subset S^N(V[x,y])$ generate a distributive lattice.
	\end{introthm}
	
	The same comparison shows that in the $y$-graded components of the common quotient of $\mathcal H_j$ and $\mathcal H_{j-1}$ the cohomology ring $H^*(\Gr(k,j))$ of Theorem~\ref{thm:introB} is replaced by $H^*(\Gr(k,j-1))$, while the module $\GPD_{j,k}$ stays the same. Summing the characters over $k$ we obtain the upper bound in the $\frac{N!}2$-conjecture; the lower bound was explained in our previous work~\cite{KM}.
	
	\begin{introthm}[Theorems~\ref{thm:N!/2} and~\ref{thm:N!/2-gl}]\label{thm:introD}
		Let $\lambda=(a,b)\vdash N+1$ with $a>b\ge1$, and let $\nu_1=(a-1,b)$, $\nu_2=(a,b-1)$. Then
		\[
		\dim\big(\GH^{\mathfrak S}_{\nu_1}\cap\GH^{\mathfrak S}_{\nu_2}\big)=\frac{N!}{2}.
		\]
		Moreover, the common quotient $S^N(V[x,y])/(\cJ_{\nu_1}+\cJ_{\nu_2})$ is isomorphic to $S^2V\otimes V^{\otimes(N-2)}$ as a $\gl_V$-module.
	\end{introthm}
	
	In Section~\ref{sec:two-rows} we state and prove Theorems~\ref{thm:introC} and~\ref{thm:introD}; the equivalence with the statements about the symmetric group is the Schur--Weyl duality of~\cite[\S3.3]{KM}.
	Theorem~\ref{thm:introC} motivates the following conjecture:
	
	\begin{conjecture}\label{conj:common-basis}
		For every $N$ there is a basis of $\mathbb C[x_1,y_1,\dots,x_N,y_N]$, consisting of bihomogeneous polynomials, such that every Garsia--Haiman module $\GH^{\mathfrak S}_\lambda$, $\lambda\vdash N$, is spanned by a subset of this basis. Equivalently, the Garsia--Haiman modules $\GH^{\mathfrak S}_\lambda$, $\lambda\vdash N$, generate a distributive lattice of subspaces of $\mathbb C[x_1,y_1,\dots,x_N,y_N]$.
	\end{conjecture}
	In this paper the conjecture is proved for the diagrams with at most two rows (Corollary~\ref{cor:common-basis}). Moreover, the conjecture is known for hook shapes after~\cite{Adin-GS-hooks},\cite{Arm} and was verified by computer for $N\le5$.
	However, we do not know an explicit common basis even for two rows diagrams.
	The conjecture is in the spirit of the Science Fiction conjecture of Bergeron and Garsia~\cite{BG}, which predicts bases of Garsia--Haiman modules compatible with the intersections of the modules $\GH^{\mathfrak S}_{\nu_i}$ for the diagrams $\nu_i$ obtained by removing corners from one diagram.
	
	\subsection{Organization of the paper}
	
	Section~\ref{sec:fusion} recalls fusion products, the basic examples and the description of Garsia--Haiman modules as iterated fusion products. Section~\ref{sec:GPD} is devoted to generalized parabolic Demazure modules: the short exact sequences, the monomial relations and the Ext-vanishing. Section~\ref{sec:two-rows} contains the presentation of the modules $\mathcal H_j^k$ and the proofs of Theorems~\ref{thm:introC} and~\ref{thm:introD}. Appendix~\ref{app:symm} collects facts on symmetric functions and plethysm. In Appendix~\ref{app:categories} we explain the plethystic substitution $X\mapsto X/(1-q)$ of Haiman's theory on the representation-theoretic side: it is the character of the functor $V\mapsto V[t]$, and the Schur functors $S_\kappa(V[t])$ are the projective covers in the category of polynomial representations of $\gl_V[t]$.
	
	\medskip\noindent\textbf{Acknowledgements.}
	We thank Boris Feigin, Ilya Dumanski and Syu Kato for stimulating conversations.
	
	{\textbf{\it{Statement on the use of AI:}}}
	We used AI tools (Claude Sonet 5.5) to proofread and polish this manuscript and to draft parts of the text following the authors' instructions. The final version of the text was proofread and edited by the authors alone.

	\section{Fusion product}\label{sec:fusion}
	
	\subsection{Definition}\label{ssec:def-fusion}
	
	The ground field is $\mathbb C$ (any field of characteristic zero would also be fine).
	
	\begin{notation}\label{not:conv}
		A partition $\lambda=(\lambda_0\ge\lambda_1\ge\dots\ge\lambda_l>0)$ is indexed from $0$; its Young diagram consists of the cells $(a,b)$, $0\le b\le l$, $0\le a<\lambda_b$, drawn in the first quadrant, so that the row $b$ has $\lambda_b$ cells and lies along the $x$-axis; the cell $(a,b)$ is attached to the monomial $x^ay^b$. We write $|\lambda|=\sum_b\lambda_b$.
		
		$V=\mathbb{C}^m$ has the basis $e_1,\dots,e_m$, $\gl=\gl(V)$, $E_{ab}$ are the matrix units, $\fh$ the diagonal matrices, $\eps_a$ the weights, $\omega_k=\eps_1+\dots+\eps_k$ ($\omega_0=0$); $\succeq$ is the dominance order. The dimension $m$ is large compared to the size of the partitions ($m\ge 3N+1$ for the results quoted from~\cite{KM}; the new arguments of Section~\ref{sec:two-rows} use only $m\ge N$). For a Lie algebra $\fg$ we write $\fg[t]=\fg\otimes\mathbb{C}[t]$, $at^k:=a\otimes t^k$, and $\U(\fg[t])_{\le n}$ for the span of the products $(a_1t^{r_1})\cdots(a_pt^{r_p})$ with $r_1+\dots+r_p\le n$.
	\end{notation}
	
	\begin{definition}[Feigin--Loktev~\cite{FeLo}]\label{def:fusion}
		Let $M_1,\dots,M_l$ be graded cyclic $\fg[t]$-modules with cyclic vectors $m_i$ (modules over $\fg$ are regarded as graded $\fg[t]$-modules with $\fg t^{k}$, $k>0$, acting by zero). For $\alpha\in\mathbb{C}$ let $M(\alpha)$ be the module $M$ with the action twisted by the automorphism $at^k\mapsto a(t+\alpha)^k$ of $\fg[t]$; for a $\fg$-module this is the evaluation module $at^k\mapsto\alpha^ka$. For pairwise distinct $\alpha_1,\dots,\alpha_l$ the tensor product $T=M_1(\alpha_1)\otimes\dots\otimes M_l(\alpha_l)$ is cyclic with the cyclic vector $m_1\otimes\dots\otimes m_l$~\cite{FeLo}. It carries the fusion filtration
		\[
		F_n(T):=\U(\fg[t])_{\le n}\cdot(m_1\otimes\dots\otimes m_l),\qquad F_{-1}=0,
		\]
		and the fusion product is the associated graded module
		\[
		M_1\star_t\dots\star_tM_l(\alpha_1,\dots,\alpha_l):=\bigoplus_{n\ge0}F_n(T)/F_{n-1}(T).
		\]
	\end{definition}
	
	The fusion product is a graded cyclic $\fg[t]$-module generated by the class of $m_1\otimes\dots\otimes m_l$; as a $\fg$-module it is isomorphic to $M_1\otimes\dots\otimes M_l$ (the subspaces $F_n$ are $\fg$-submodules), so $\dim=\prod\dim M_i$. When the isomorphism class does not depend on $\alpha_1,\dots,\alpha_l$ we omit them. It is not known in general whether the fusion product depends on the parameters, nor whether it is associative; both were conjectured to hold in~\cite{FeLo}. The cases where they are known are listed in Section~\ref{ssec:examples}.
	
	\begin{remark}[fusion over {$\fg[x]$}]\label{rem:fusion-xy}
		If $M_i$ are graded cyclic modules over the current algebra $\fg[x]$, we may regard them as modules over the Lie algebra $\fg[x]$ and fuse them with respect to a new variable $y$; the result is a module over $\fg[x][y]=\fg[x,y]$, written $M_1\star_yM_2\star_y\cdots$. The filtration is by the $y$-degree; if the $M_i$ are graded by the $x$-degree and the $m_i$ are homogeneous, then the $F_n$ are $x$-graded (the operators $(a x^ry^s)$ act on $T$ as $\sum_i\alpha_i^s(ax^r)^{(i)}$, hence are homogeneous of $x$-degree $r$). Thus the fusion product is bigraded, and its component of $y$-degree $k$ is a graded $\fg[x]$-module. This is the situation of Sections~\ref{ssec:GH-fusion} and~\ref{sec:two-rows}.
	\end{remark}
	
	\subsection{Basic known examples}\label{ssec:examples}
	
	\subsubsection{Local Weyl modules}\label{sssec:weyl}
	Let $\fg$ be a simple Lie algebra (or $\gl_m$), $\lambda$ a dominant weight. The global Weyl module $\bW_\lambda$ is the $\fg[x]$-module generated by $w_\lambda$ with the defining relations
	\[
	\fn_+[x]\,w_\lambda=0,\qquad hw_\lambda=\lambda(h)w_\lambda\ (h\in\fh),\qquad (e_{-\alpha_i})^{\langle\lambda,\alpha_i^\vee\rangle+1}w_\lambda=0,
	\]
	(and, for $\gl_m$, $(\mathrm{Id}\otimes x^t)w_\lambda=0$ for $t\ge1$, i.e.\ the central currents act by zero). The local Weyl module $W_\lambda$ is its quotient by the additional relations $(h\otimes x^t)w_\lambda=0$, $h\in\fh$, $t\ge1$. Both are graded by the $x$-degree. More details are collected in Appendix~\ref{app:categories}.
	
	\begin{theorem}[Chari--Loktev~\cite{CL}, Fourier--Littelmann~\cite{FL1,FL2}, Naoi~\cite{N}]\label{thm:Weyl-fusion}
		Let $\lambda=\sum_i m_i\omega_i$ be dominant. Then, for every choice of pairwise distinct parameters,
		\[
		W_\lambda\ \simeq\ \underbrace{V_{\omega_1}\star_x\dots\star_xV_{\omega_1}}_{m_1}\star_x\dots\star_x\underbrace{V_{\omega_r}\star_x\dots\star_xV_{\omega_r}}_{m_r},
		\]
		in particular the fusion product does not depend on the parameters. For $\gl_V$ and $\lambda=r\eps_1$ one gets $W_{r\eps_1}\simeq V\star_x\dots\star_xV$ ($r$ factors), a module of dimension $m^r$.
	\end{theorem}
	
	For simply laced $\fg$ the independence of the parameters is explained by the coincidence of local Weyl modules with the Demazure modules of level one: $W_\lambda\simeq D(1,\lambda)$~\cite{FL2}, and the fusion product of Demazure modules of the same level is again a Demazure module, $D(\ell,\ell\lambda)\star D(\ell,\ell\lambda')\simeq D(\ell,\ell(\lambda+\lambda'))$~\cite{FL1}; the right-hand side does not depend on the parameters.
	
	\subsubsection{Fusion products of Demazure modules of different levels}\label{sssec:dumanski}
	For a simply laced $\fg$, a dominant weight $\lambda$ and $\ell\ge1$ let $D(\ell,\lambda)$ be the $\fg[x]$-stable affine Demazure module of level $\ell$ (see~\cite[\S1.3]{D}). Recall that $D(1,\lambda)\simeq W_\lambda$, that $D(\ell,\ell\omega_k)$ is the Kirillov--Reshetikhin module, and that for $\fg=\mathfrak{sl}_n$ one has $D(\ell,\ell\omega_k)\simeq V_{\ell\omega_k}$ as a $\fg$-module.
	
	\begin{theorem}[Dumanski~{\cite[Cor.~2.3]{D}}]\label{thm:Dumanski}
		Let $\fg$ be simply laced, $\lambda_1,\dots,\lambda_p$ and $\mu$ dominant weights, $\ell_1\ge\dots\ge\ell_p\ge\ell\ge0$. Then the fusion product
		\[
		D(\ell_1,\ell_1\lambda_1)\star_x\dots\star_xD(\ell_p,\ell_p\lambda_p)\star_xD(\ell,\mu)
		\]
		does not depend on the parameters, is associative, and is isomorphic to a generalized Demazure module of Naoi~\cite{Naoi} (a generalized parabolic Demazure module in the terminology of Section~\ref{sec:GPD}). Equivalently, it is the cyclic submodule generated by the tensor product of the cyclic vectors in
		\[
		D(\ell_1-\ell_2,(\ell_1-\ell_2)\lambda_1)\otimes D(\ell_2-\ell_3,(\ell_2-\ell_3)(\lambda_1+\lambda_2))\otimes\dots\otimes D(\ell,\ell(\lambda_1+\dots+\lambda_p)+\mu).
		\]
	\end{theorem}
	
	The independence of the parameters comes from the geometry of the Beilinson--Drinfeld Grassmannian: the dual of the space of sections of a line bundle on the global convolution diagram (a flat family over $\mathbb A^p$) is, over the diagonal, the left-hand side with parameters collapsed, and over generic points the tensor product of the factors; the fusion filtration is the specialization (see~\cite[Thm.~2.1]{D}). Taking $\ell_i-\ell_{i+1}=1$, $\ell=0$, $\mu=0$ one obtains exactly the generalized parabolic Demazure modules of Section~\ref{sec:GPD}:
	\[
	D(p,p\lambda_1)\star D(p-1,(p-1)\lambda_2)\star\dots\star D(1,\lambda_p)\ \simeq\ \GPD(\lambda_1,\dots,\lambda_p)\subset W_{\lambda_1}\otimes W_{\lambda_1+\lambda_2}\otimes\dots\otimes W_{\lambda_1+\dots+\lambda_p}.
	\]
	The details are given in Section~\ref{sec:GPD}.
	
	\subsection{Garsia--Haiman modules as iterated fusion products}\label{ssec:GH-fusion}
	
	Let $\lambda\vdash N$ be a Young diagram with cells $(a_c,b_c)$, $c=1,\dots,N$, and
	\[
	\Delta_\lambda=\det\big(x_i^{a_c}y_i^{b_c}\big)_{i,c=1}^N\in\mathbb{C}[x_1,y_1,\dots,x_N,y_N].
	\]
	The Garsia--Haiman module $\GH^{\mathfrak S}_\lambda$ is the span of all partial derivatives of $\Delta_\lambda$, and $J_\lambda=\{f:\ f(\partial)\Delta_\lambda=0\}$ is its annihilator ideal~\cite{GH}. By Haiman's theorem~\cite{Haiman} (the $N!$-conjecture) $\dim\GH^{\mathfrak S}_\lambda=N!$.
	
	For $m\ge N$ there is a version of Schur--Weyl duality for the current algebra $\gl_V[x,y]$~\cite[\S3]{KM}: the functor of the $(1^N)$-weight space identifies the category of polynomial $\gl_V[x,y]$-modules of degree $N$ with the category of modules over $\mathfrak S_N\ltimes\mathbb{C}[x_1,y_1,\dots,x_N,y_N]$; the module $S^N(V[x,y])$ ($V[x,y]=V\otimes\mathbb{C}[x,y]$) corresponds to the polynomial ring, and $\mathfrak S_N$-stable ideals correspond to $\gl_V[x,y]$-submodules. Let $\cJ_\lambda\subset S^N(V[x,y])$ be the submodule corresponding to $J_\lambda$, and
	\[
	\GH^{\gl}_\lambda:=S^N(V[x,y])/\cJ_\lambda .
	\]
	It is a bigraded cyclic $\gl_V[x,y]$-module with the cyclic vector $e_1^N$; its $(1^N)$-weight space is $\mathbb{C}[\bx,\by]/J_\lambda$. 
	
	\begin{theorem}[{\cite[Theorem A]{KM}}]\label{thm:GH-fusion}
		For every Young diagram $\lambda\vdash N$ there is a canonical isomorphism of bigraded $\gl_V[x,y]$-modules 
		$$\mathcal H_\lambda:= 
		\mathcal H_\lambda:=W_{\lambda_0\eps_1}(x)\star_yW_{\lambda_1\eps_1}(x)\star_y\cdots\star_yW_{\lambda_l\eps_1}(x)
		\twoheadrightarrow\GH^{\gl}_\lambda$$
		from the iterated fusion product of $N$ copies of $V$ to the Garsia-Haiman module.
		Which sends the cyclic vector to the class of $e_1^N$.
	\end{theorem}
	
	\begin{remark}
		In particular, the independence of the fusion parameters for $\mathcal H_\lambda$ follows from Haiman's $N!$-theorem; we do not know a proof avoiding the geometry of the isospectral Hilbert scheme. In the two-row case $\lambda=(\lambda_0,\lambda_1)$ we obtain
		\[
		\GH^{\gl}_{(\lambda_0,\lambda_1)}\simeq W_{\lambda_0\eps_1}(x)\star_yW_{\lambda_1\eps_1}(x),
		\]
		and the $\gl_V$-module structure is $V^{\otimes N}$, $\dim=m^N$. Sections~\ref{sec:GPD} and~\ref{sec:two-rows} describe this module in detail.
	\end{remark}
	
	\section{Generalized parabolic Demazure modules and the recursion for them}\label{sec:GPD}
	
	In this section $\fg=\gl_m$ (everything below holds for any simply laced $\fg$; polynomial $\gl_m$-modules are regarded as $\mathfrak{sl}_m$-modules with the central currents acting by zero). We write $\mu=\sum_k\mu_k\omega_k$ for dominant weights and $f_k=e_{-\alpha_k}=E_{k+1,k}$. The modules $W_\mu$ are the local Weyl modules, with the cyclic (highest weight) vectors $w_\mu$.
	
	\subsection{Generalized parabolic Demazure modules}\label{ssec:GPD-def}
	
	\begin{definition}\label{def:GPD}
		Let $\mu^1,\dots,\mu^p$ be dominant weights. The generalized parabolic Demazure module
		\[
		\GPD(\mu^1,\dots,\mu^p)\ \subset\ W_{\mu^1}\otimes W_{\mu^1+\mu^2}\otimes\dots\otimes W_{\mu^1+\dots+\mu^p}
		\]
		is the $\fg[x]$-submodule generated by the tensor product $w_{\mu^1}\otimes w_{\mu^1+\mu^2}\otimes\dots\otimes w_{\mu^1+\dots+\mu^p}$ of the cyclic vectors.
	\end{definition}
	
	It is a graded cyclic $\fg[x]$-module, generated by a vector of weight $\sum_{i=1}^p(p-i+1)\mu^i$. For $p=1$ it is the local Weyl module $W_{\mu^1}$, and $\GPD(0,\nu)=W_\nu$.
	
	\begin{remark}\label{rem:GPD-names}
		In Naoi's terminology~\cite{Naoi} a generalized Demazure module is the $\hat{\fb}$-submodule $D(\xi_1,\dots,\xi_p)=\U(\hat\fb)(v_{\xi_1}\otimes\dots\otimes v_{\xi_p})$ of a tensor product of integrable highest weight modules of the affine Lie algebra, generated by the product of extremal vectors. The module $\GPD(\mu^1,\dots,\mu^p)$ is the generalized Demazure module $D(w_0\mu^1,w_0(\mu^1+\mu^2),\dots,w_0(\mu^1+\dots+\mu^p))$ with anti-dominant extremal weights of level one (note that $W_\lambda\simeq D(1,\lambda)$ is a level one Demazure module~\cite{FL2}). Such modules are stable under the whole parabolic subalgebra $\fg[x]$ (a generalized Demazure module attached to arbitrary extremal weights is in general only a $\hat\fb$-module), whence the name generalized parabolic Demazure modules. We use the normalization in which all cyclic vectors are highest weight vectors; the passage from the extremal vectors of weight $w_0\lambda$ is a twist by an automorphism of $\fg[x]$.
	\end{remark}
	
	\begin{theorem}[Dumanski]\label{thm:GPD-fusion}
		For dominant $\mu^1,\dots,\mu^p$ we have an isomorphism of graded $\fg[x]$-modules
		\[
		\GPD(\mu^1,\dots,\mu^p)\ \simeq\ D(p,p\mu^1)\star_xD(p-1,(p-1)\mu^2)\star_x\dots\star_xD(1,\mu^p),
		\]
		sending the cyclic vector to the cyclic vector; the fusion product does not depend on the parameters and is associative. Consequently, as $\gl_m$-modules (writing $\mu^i=\sum_k\mu^i_k\omega_k$)
		\[
		\GPD(\mu^1,\dots,\mu^p)\ \simeq\ \bigotimes_{i=1}^p\bigotimes_{k}\big(V_{(p-i+1)\omega_k}\big)^{\otimes\mu^i_k},\qquad
		\ch\GPD(\mu^1,\dots,\mu^p)=\prod_{i=1}^p\prod_k\big(\ch V_{(p-i+1)\omega_k}\big)^{\mu^i_k}.
		\]
		Here $V_{\ell\omega_k}$ is the irreducible module whose Young diagram is the rectangle with $k$ rows and $\ell$ columns.
	\end{theorem}
	\begin{proof}
		The first statement is Theorem~\ref{thm:Dumanski} for $\ell_i=p-i+1$, $\ell=0$, $\mu=0$, in the description as a cyclic submodule of the tensor product of level one Demazure modules $D(1,\mu^1+\dots+\mu^i)\simeq W_{\mu^1+\dots+\mu^i}$. The fusion product is isomorphic to the tensor product as a $\fg$-module, $D(\ell,\ell\mu)\simeq\bigotimes_kD(\ell,\ell\omega_k)^{\otimes\mu_k}$ as $\fg$-modules (the fusion of Demazure modules of the same level is a Demazure module~\cite{FL1}), and $D(\ell,\ell\omega_k)\simeq V_{\ell\omega_k}$ in type $A$.
	\end{proof}
	
	\begin{example}\label{ex:GPD-two-rows}
		For a two-row partition $\lambda=(\lambda_0\ge\lambda_1)$, $N=|\lambda|$, and $0\le k\le\lambda_1$ put
		\[
		\GPD_{\lambda,k}:=\GPD\big((\lambda_1-k)\omega_1,\ (\lambda_0-\lambda_1)\omega_1+k\omega_2\big)\subset W_{(\lambda_1-k)\eps_1}\otimes W_{\lambda_0\eps_1+k\eps_2}.
		\]
		Its cyclic vector has weight $(N-k)\eps_1+k\eps_2$, $\GPD_{\lambda,\lambda_1}=W_{\lambda_0\eps_1+\lambda_1\eps_2}$, and by Theorem~\ref{thm:GPD-fusion}
		\[
		\ch\GPD_{\lambda,k}=h_2^{\lambda_1-k}\,p_1^{\lambda_0-\lambda_1}\,e_2^{k},\qquad \dim\GPD_{\lambda,k}=m^{\lambda_0}\Big(\frac{m+1}2\Big)^{\lambda_1-k}\Big(\frac{m-1}2\Big)^{k},
		\]
		because $V_{2\omega_1}=S^2V$, $V_{\omega_1}=V$, $V_{\omega_2}=\Lambda^2V$. In particular $\sum_{k=0}^{\lambda_1}\binom{\lambda_1}k\dim\GPD_{\lambda,k}=m^{\lambda_0}\,m^{\lambda_1}=m^N$.
		When $N$ is fixed and $\lambda=(N-j,j)$ we write $\GPD_{j,k}:=\GPD_{(N-j,j),k}$; thus $\ch\GPD_{j,k}=h_2^{j-k}p_1^{N-2j}e_2^k$.
	\end{example}

	Not much is known about generalized parabolic Demazure modules. Their dimensions, and even their $\gl_m$-characters, are known: by Theorem~\ref{thm:GPD-fusion} these modules are fusion products of irreducible modules, and the fusion product does not change the $\fg$-module structure. Their graded characters, i.e.\ the characters that also record the $x$-degree, are not known in general. For $p=1$ the graded character of $\GPD(\mu)=W_\mu$ is the specialized Macdonald polynomial $P_\mu(X;q,0)$. For $p=2$ we find the defining relations of $\GPD(\mu,\nu)$ (Corollary~\ref{cor:monomial}) and show that $\GPD(\mu,\nu)$ coincides with the module $M(\nu,2\mu)$ of Wand and Biswal--Chari--Shereen--Wand~\cite{Wand,BCSW}; for these modules both the defining relations and the graded characters were already known (Remark~\ref{rem:BCSW}). For $p\ge3$ we know neither a formula for the graded character nor a complete list of defining relations.
	
	Our proof of the presentation is short and uses fusion products; it also gives the exact sequences~\eqref{eq:SES} for arbitrary $p$. The family $\GPD_{\lambda,k}$ of Example~\ref{ex:GPD-two-rows} is the one needed in Section~\ref{sec:two-rows}: the $y$-graded components of the two-row Garsia--Haiman modules have filtrations whose subquotients are isomorphic to these modules (Corollary~\ref{cor:A}).
	
	\subsection{Short exact sequences}\label{ssec:SES}

	\begin{theorem}\label{thm:SES}
		Let $\mu^1,\dots,\mu^p$ be dominant weights, $p\ge2$, and $1\le k<m$ (we put $\omega_0=0$). Then there is a short exact sequence of graded $\fg[x]$-modules
		\begin{equation}\label{eq:SES}
			\begin{aligned}
				0\to{}&\GPD(\mu^1,\dots,\mu^{p-1},\mu^p+\omega_{k-1}+\omega_{k+1})\xrightarrow{\ \iota\ }\GPD(\mu^1,\dots,\mu^{p-1},\mu^p+2\omega_k)\\
				&\xrightarrow{\ \pi\ }\GPD(\mu^1,\dots,\mu^{p-2},\mu^{p-1}+\omega_k,\mu^p)\to0.
			\end{aligned}
		\end{equation}
		Here $\pi$ sends the cyclic vector to the cyclic vector, and $\iota$ sends the cyclic vector to $f_kx^{s}$ applied to the cyclic vector of the middle module, where $s=\langle\mu^1+\dots+\mu^p,\alpha_k^\vee\rangle+1$.
	\end{theorem}
	
	\begin{proof}
		Put $\Sigma=\mu^1+\dots+\mu^p$, $\lambda=\Sigma+2\omega_k$, $\lambda'=\lambda-\alpha_k=\Sigma+\omega_{k-1}+\omega_{k+1}$ and $s=\langle\lambda,\alpha_k^\vee\rangle-1=\langle\Sigma,\alpha_k^\vee\rangle+1$.
		
		Since $\lambda'\preceq\lambda$, we have $D(1,\lambda')\subset D(1,\lambda)$, i.e.\ an injective homomorphism of $\fg[x]$-modules $\iota_0\colon W_{\lambda'}\to W_\lambda$ (via $W_\lambda\simeq D(1,\lambda)$~\cite{FL2}). The image of the cyclic vector is a vector $v'\in W_\lambda$ of weight $\lambda'$, homogeneous of degree equal to the difference of the energies of the extremal vectors, $\tfrac12\big((\lambda,\lambda)-(\lambda',\lambda')\big)=\langle\lambda,\alpha_k^\vee\rangle-1=s$ (we used $(\alpha_k,\alpha_k)=2$). The weight space $W_\lambda[\lambda-\alpha_k]$ is spanned by the vectors $f_kx^tw_\lambda$, so $v'=c\,f_kx^sw_\lambda$ with $c\ne0$ because $\iota_0$ is injective. Thus
		\[
		\iota_0(w_{\lambda'})=c\,f_kx^{s}w_\lambda,\qquad c\neq0.
		\]
		
		Let $\iota=\id\otimes\dots\otimes\id\otimes\iota_0$ on the ambient tensor products; it is injective. The cyclic vector of the left module of~\eqref{eq:SES} goes to $w_{\mu^1}\otimes\dots\otimes w_{\mu^1+\dots+\mu^{p-1}}\otimes c\,f_kx^sw_\lambda$. By the local Weyl relations $f_kx^tw_\nu=0$ for $t\ge\langle\nu,\alpha_k^\vee\rangle$, and $s\ge\langle\mu^1+\dots+\mu^i,\alpha_k^\vee\rangle$ for $i<p$, therefore by the Leibniz rule this vector equals $c\,f_kx^s$ applied to the cyclic vector of the middle module. Hence $\iota$ maps the left module into the middle one, injectively, and its image is $\U(\fg[x])f_kx^sv$, where $v$ is the cyclic vector of the middle module.
		
		The module $\GPD(\omega_k,\nu)$ is generated by the vector $w_{\omega_k}\otimes w_{\omega_k+\nu}$ of weight $\nu+2\omega_k$, annihilated by $\fn_+[x]$ and $\fh\otimes x\mathbb{C}[x]$, and it is integrable; hence there is a surjection $\varphi\colon W_{\nu+2\omega_k}\twoheadrightarrow\GPD(\omega_k,\nu)$ sending cyclic vector to cyclic vector. By Theorem~\ref{thm:GPD-fusion} (independence of parameters and associativity),
		\[
		\GPD(\mu^1,\dots,\mu^{p-1},\mu^p+2\omega_k)\simeq D(p,p\mu^1)\star\dots\star D(2,2\mu^{p-1})\star W_{\mu^p+2\omega_k},
		\]
		\[
		\GPD(\mu^1,\dots,\mu^{p-1}+\omega_k,\mu^p)\simeq D(p,p\mu^1)\star\dots\star D(2,2\mu^{p-1})\star\big[D(2,2\omega_k)\star D(1,\mu^p)\big],
		\]
		where for the second isomorphism we used $D(2,2\mu^{p-1})\star D(2,2\omega_k)\simeq D(2,2(\mu^{p-1}+\omega_k))$ (fusion of Demazure modules of the same level~\cite{FL1}) and $\GPD(\omega_k,\mu^p)=D(2,2\omega_k)\star D(1,\mu^p)$. Since the fusion product preserves surjections of cyclic modules~\cite[Proposition ``Fusion product preserves surjections'']{KM}, $\id\star\dots\star\id\star\varphi$ gives a surjection $\pi$, which maps the cyclic vector to the cyclic vector.
		
		Notice that $\pi\iota=0$. Indeed $\pi(f_kx^sv)=f_kx^sv'$, where $v'=w_{\mu^1}\otimes\dots\otimes w_{\mu^1+\dots+\mu^{p-2}}\otimes w_{\mu^1+\dots+\mu^{p-1}+\omega_k}\otimes w_{\Sigma+\omega_k}$ is the cyclic vector of the right module, and $f_kx^s$ kills each tensor factor, because $s\ge\langle\text{weight},\alpha_k^\vee\rangle$ for each of them (the last one has $\langle\Sigma+\omega_k,\alpha_k^\vee\rangle=s$). Finally, by Theorem~\ref{thm:GPD-fusion} and the character equality $\ch V_{\omega_k}\cdot\ch V_{\omega_k}=\ch V_{2\omega_k}+\ch V_{\omega_{k-1}}\cdot\ch V_{\omega_{k+1}}$, the characters of the three modules add up: the middle one contains the factor $(\ch V_{\omega_k})^2$ in place of $\ch V_{\omega_{k-1}}\ch V_{\omega_{k+1}}$ (left) and $\ch V_{2\omega_k}$ (right, at level two), all other factors being the same. Hence $\dim\ker\pi=\dim(\text{left module})=\dim\operatorname{im}\iota$, and $\operatorname{im}\iota\subset\ker\pi$ forces equality.
	\end{proof}
	\begin{example}
		For $p=2$, $\mu^1=\mu$, $\mu^2=\nu$ the theorem reads as
		\[
		0\to\GPD(\mu,\nu+\omega_{k-1}+\omega_{k+1})\to\GPD(\mu,\nu+2\omega_k)\to\GPD(\mu+\omega_k,\nu)\to0,
		\]
		and for $\mu=0$ it reads $0\to W_{\nu+\omega_{k-1}+\omega_{k+1}}\to W_{\nu+2\omega_k}\to\GPD(\omega_k,\nu)\to0$.
		
	\end{example}

	\begin{corollary}\label{cor:monomial}
		Under the assumptions of Theorem~\ref{thm:SES}, with $s=\langle\mu^1+\dots+\mu^p,\alpha_k^\vee\rangle+1$,
		\[
		\GPD(\mu^1,\dots,\mu^{p-1}+\omega_k,\mu^p)\ \simeq\ \GPD(\mu^1,\dots,\mu^{p-1},\mu^p+2\omega_k)\Big/\U(\fg[x])\,f_kx^{s}\,v ,
		\]
		where $v$ is the cyclic vector of $\GPD(\mu^1,\dots,\mu^{p-1},\mu^p+2\omega_k)$. In particular for $p=2$ and dominant $\mu=\sum_i m_i\omega_i$, $\nu$:
		\begin{equation}\label{eq:monomial}
			\GPD(\mu,\nu)\ \simeq\ W_{2\mu+\nu}\Big/\U(\fg[x])\Big\{\,f_ix^{\langle\mu+\nu,\alpha_i^\vee\rangle}\,w_{2\mu+\nu}\ \Big|\ m_i>0\Big\}.
		\end{equation}
		Thus $\GPD(\mu,\nu)$ is the cyclic module generated by $w$ with the defining relations of the local Weyl module $W_{2\mu+\nu}$ and with the additional monomial relations $f_ix^{n_i}w=0$, $n_i=\langle\mu+\nu,\alpha_i^\vee\rangle$, i.e.\ the vanishing on $w$ of single root vectors $e_{-\alpha_i+n_i\delta}$ of the affine Lie algebra, for the simple roots with $m_i>0$.
	\end{corollary}
	\begin{proof}
		The first statement is the exactness of~\eqref{eq:SES} together with Step 2. For~\eqref{eq:monomial} we use induction on $|\mu|=\sum m_i$: for $\mu=0$ it is $\GPD(0,\nu)=W_\nu$. Let $\mu'=\mu+\omega_k$. By the first statement and the induction hypothesis applied to $\GPD(\mu,\nu+2\omega_k)\simeq W_{2\mu+\nu+2\omega_k}/\langle f_ix^{\langle\mu+\nu,\alpha_i^\vee\rangle+2\delta_{ik}}w:m_i>0\rangle$ we get $\GPD(\mu',\nu)\simeq W_{2\mu'+\nu}$ modulo the relations $f_ix^{\langle\mu+\nu,\alpha_i^\vee\rangle}w$ ($i\ne k$, $m_i>0$), $f_kx^{\langle\mu+\nu,\alpha_k^\vee\rangle+2}w$ (if $m_k>0$) and the new relation $f_kx^{s}w$, $s=\langle\mu+\nu,\alpha_k^\vee\rangle+1=\langle\mu'+\nu,\alpha_k^\vee\rangle$. Since $\langle\mu',\alpha_i^\vee\rangle=\langle\mu,\alpha_i^\vee\rangle$ for $i\ne k$, and $f_kx^{t+1}w\in\U(\fg[x])f_kx^tw$ in every local Weyl module ($h_kx\,w=0$ and $[h_kx,f_kx^t]=-2f_kx^{t+1}$), the claim follows.
	\end{proof}

	\begin{remark}\label{rem:BCSW}
		For $\fg=\mathfrak{sl}_{n+1}$ and dominant weights $\nu,\lambda$ write $\lambda=2\lambda_0+\lambda_1$ with $\lambda_0$ dominant and $\langle\lambda_1,\alpha_i^\vee\rangle\le1$ for all $i$. J.~Wand~\cite{Wand} introduced the graded $\fg[x]$-module $M(\nu,\lambda)$ generated by a vector $w$ subject to the defining relations of the local Weyl module $W_{\nu+\lambda}$ and the relations
		\[
		e_{-\alpha}\,x^{\langle\nu,\alpha^\vee\rangle+\lceil\langle\lambda,\alpha^\vee\rangle/2\rceil}\,w=0\qquad\text{for all positive roots }\alpha
		\]
		(see~\cite[\S2.2]{BCSW} and~\cite{BCM}). These modules interpolate between the local Weyl modules $M(\nu,0)=W_\nu$ and the level two Demazure modules $M(0,\lambda)\simeq D(2,\lambda)$~\cite{CSVW}. For $\lambda=2\mu$ the exponent equals $\langle\mu+\nu,\alpha^\vee\rangle$, and we have
		\begin{equation}
			\label{eq::GPD=M}	
			\GPD(\mu,\nu)\ \simeq\ M(\nu,2\mu).
		\end{equation}
		Indeed, by~\eqref{eq:monomial} the module $M(\nu,2\mu)$ is a quotient of $\GPD(\mu,\nu)$: for simple roots its extra relations are exactly the monomial relations. Conversely, the relations of $M(\nu,2\mu)$ hold in $\GPD(\mu,\nu)\subset W_\mu\otimes W_{\mu+\nu}$ by the Leibniz rule, because $e_{-\alpha}x^tw_\kappa=0$ in $W_\kappa$ for every positive root $\alpha$ and $t\ge\langle\kappa,\alpha^\vee\rangle$~\cite{CP}, and $\langle\mu,\alpha^\vee\rangle\le\langle\mu+\nu,\alpha^\vee\rangle$. Both surjections send the cyclic vector to the cyclic vector, so they are isomorphisms. 
		In particular, we observe that in $\GPD(\mu,\nu)$ the relations for non-simple roots follow from those for simple roots.
		
		Moreover, under Identification~\eqref{eq::GPD=M} we see that a particular case $p=2$ of our the exact sequence~\eqref{eq:SES}  coincides with the exact sequence of~\cite[Proposition~2.5(i)(a)]{BCSW} for $\lambda=2\mu$.
	\end{remark}
	
	\begin{corollary}\label{cor:two-rows-presentation}
		For a two-row partition $\lambda=(\lambda_0\ge\lambda_1)$ and $0\le k\le\lambda_1$,
		\[
		\GPD_{\lambda,k}\ \simeq\ W_{(N-k)\eps_1+k\eps_2}\Big/\U(\gl_V[x])\,E_{21}x^{\lambda_0-k}\,w .
		\]
		Indeed this is~\eqref{eq:monomial} with $\mu=(\lambda_1-k)\omega_1$, $\nu=(\lambda_0-\lambda_1)\omega_1+k\omega_2$, $\langle\mu+\nu,\alpha_1^\vee\rangle=\lambda_0-k$ (for $k=\lambda_1$ the relation holds in $W_\lambda$ already).
	\end{corollary}
	
	\subsection{Ext-vanishing}\label{ssec:Ext}
	
	Recall (Appendix~\ref{app:categories}) that $\Ext^\bullet$ denotes the Ext groups in the category of graded finite dimensional $\gl_V[x]$-modules (the highest weight category of~\cite{BBCKL}), that $V_\kappa$ are the simple modules (concentrated in degree $0$), and that the local Weyl modules satisfy $\Ext^\bullet(W_\lambda,V_\kappa)=0$ for $\kappa\not\succeq\lambda$ (Lemma~\ref{lem:Ext-Weyl}).
	
	\begin{proposition}\label{prop:Ext}
		Let $\lambda=(\lambda_0\ge\lambda_1)$, $0\le k\le\lambda_1$, and let $\kappa\vdash N$ be a partition with $\kappa\not\succeq\lambda$ (i.e.\ $\kappa$ has at least three rows or $\kappa_0<\lambda_0$). Then
		\[
		\Ext^\bullet(\GPD_{\lambda,k},V_\kappa)=0 .
		\]
	\end{proposition}
	\begin{proof}
		Induction on $d=\lambda_1-k$ (for all $\lambda$ simultaneously). If $d=0$, then $\GPD_{\lambda,\lambda_1}=W_\lambda$ and the claim is Lemma~\ref{lem:Ext-Weyl}. Let $d\ge1$ and put $\lambda^+=(\lambda_0+1,\lambda_1-1)$. Theorem~\ref{thm:SES} for $p=2$, $k=1$ (the simple root $\alpha_1$, $\omega_0=0$), $\mu=(\lambda_1-k-1)\omega_1$, $\nu=(\lambda_0-\lambda_1)\omega_1+k\omega_2$ gives the exact sequence
		\[
		0\to\GPD_{\lambda,k+1}\to\GPD_{\lambda^+,k}\to\GPD_{\lambda,k}\to0
		\]
		(check: the left term is $\GPD(\mu,\nu+\omega_2)$, the middle one $\GPD(\mu,\nu+2\omega_1)$, the right one $\GPD(\mu+\omega_1,\nu)$). The modules on the left and in the middle have defect $\lambda_1-(k+1)=d-1$, resp.\ $(\lambda_1-1)-k=d-1$. If $\kappa\not\succeq\lambda$ then $\kappa\not\succeq\lambda^+$ as well (because $\lambda^+\succeq\lambda$), so by induction $$\Ext^\bullet(\GPD_{\lambda,k+1},V_\kappa)=\Ext^\bullet(\GPD_{\lambda^+,k},V_\kappa)=0.$$ 
		The long exact sequence
		\[
		\dots\to\Ext^i(\GPD_{\lambda,k},V_\kappa)\to\Ext^i(\GPD_{\lambda^+,k},V_\kappa)\to\Ext^i(\GPD_{\lambda,k+1},V_\kappa)\to\Ext^{i+1}(\GPD_{\lambda,k},V_\kappa)\to\dots
		\]
		gives $\Ext^\bullet(\GPD_{\lambda,k},V_\kappa)=0$.
	\end{proof}
	
	\begin{corollary}\label{cor:char-span}
		The graded character of $\GPD_{\lambda,k}$ lies in $\mathbb Z[q]\langle\,s_\kappa\big(\tfrac{X}{1-q}\big):\ \kappa\succeq\lambda\,\rangle$ (Appendix~\ref{app:categories}, Corollary~\ref{cor:Ext-char}).
	\end{corollary}
	
	\section{Garsia--Haiman modules with two rows}\label{sec:two-rows}
	
	\subsection{Overview}\label{ssec:overview}
	
	Let $\lambda=(N-j,j)$ be a two-row partition of $N$ (so $0\le j\le\lfloor N/2\rfloor$). By Theorem~\ref{thm:GH-fusion} the Garsia--Haiman module of $\lambda$ is the $y$-fusion product of two local Weyl modules of $\gl_V[x]$:
	\[
	\mathcal H_j:=\GH^{\gl}_{(N-j,j)}\ \simeq\ W_{(N-j)\eps_1}(x)\star_yW_{j\eps_1}(x),\qquad \dim\mathcal H_j=m^N .
	\]
	It is a bigraded $\gl_V[x,y]$-module (where $y$-degree is the fusion degree), so
	$\mathcal H_j=\bigoplus_{k=0}^{j}\mathcal H_j^k$, where each  $\mathcal H_j^k$ is a graded $\gl_V[x]$-module. This section describes these components. The answer, in a nutshell:
	\begin{itemize}
		\item $\mathcal H_j^k$ is a quotient of the global Weyl module $\bW_{(N-k)\eps_1+k\eps_2}$ by two groups of relations: relations of the Cartan currents, which cut the algebra generated by the Cartan currents acting on the cyclic vector down to the ring
		\[
		\mathbb C[t^1_1,\dots,t^1_{j-k},t^2_1,\dots,t^2_k]^{\sS_{j-k}\times\sS_k}\big/\big(p_r(t^1_\bullet)+p_r(t^2_\bullet):\ r\ge1\big)\ \simeq\ H^*(\Gr(k,j);\mathbb C).
		\]
		and a Cartan-twisted relation of the generalized parabolic Demazure module $\GPD_{j,k}$ (Theorem~\ref{thm:presentation}). 
		Here $p_r(t^1_\bullet)$ denotes the Newton's power sums in variables $\{t^1_l\colon l=1\ldots j-k\}$ (resp. $p_r(t^1_\bullet)$ are Newton's sum for the remaining $t^2_l$-variables)
		
		\item Consequently $\mathcal H_j^k$ has a filtration by $\binom jk$ subquotients, each isomorphic to the generalized parabolic Demazure module $\GPD_{j,k}$ of Example~\ref{ex:GPD-two-rows} (Corollary~\ref{cor:A}).
	\end{itemize}
	This description allows us to compare the modules $\mathcal H_j$ for different $j$ (with $N$ fixed), and implies two main results.
	
	\begin{theorem}
		\label{thm:lattice}
		Let $N\ge2$ and $0\le j\le\lfloor N/2\rfloor$. The Garsia--Haiman ideals $J_{(N-j,j)}\subset\mathbb{C}[x_1,y_1,\dots,x_N,y_N]$ of the two-row diagrams of size $N$ generate a distributive lattice of ideals. Equivalently, the spaces $\GH^{\mathfrak S}_{(N-j,j)}$ generate a distributive lattice of subspaces of $\mathbb{C}[\bx,\by]$.	
	\end{theorem}
	
	\begin{theorem}\label{thm:N!/2}
		Let $\lambda\vdash N+1$ be a two-row diagram with two corners, i.e.\ $\lambda=(a,b)$, $a>b\ge1$, and let $\nu_1=(a-1,b)$, $\nu_2=(a,b-1)$ be the diagrams obtained by removing the corners. Then
		\[
		\dim\big(\GH^{\mathfrak S}_{\nu_1}\cap\GH^{\mathfrak S}_{\nu_2}\big)=\frac{N!}2 .
		\]
	\end{theorem}
	
	Theorems~\ref{thm:lattice} and~\ref{thm:N!/2} are stated in their original form, in terms of the symmetric group. Below we formulate and prove their Schur--Weyl dual forms, Theorems~\ref{thm:lattice-gl} and~\ref{thm:N!/2-gl}, which concern the $\gl_V[x,y]$-submodules $\cJ_{(N-j,j)}\subset S^N(V[x,y])$ and their common quotients; the symmetric group does not appear in the proofs. The equivalence of the two forms is the Schur--Weyl duality for current algebras~\cite[\S3.3]{KM}: the functor of the $(1^N)$-weight space identifies the lattice of $\gl_V[x,y]$-submodules of $S^N(V[x,y])$ with the lattice of $\sS_N$-stable ideals of $\mathbb C[\bx,\by]$, $\cJ_\nu\mapsto J_\nu$, and the common quotient $S^N(V[x,y])/(\cJ_{\nu_1}+\cJ_{\nu_2})$ with $\mathbb C[\bx,\by]/(J_{\nu_1}+J_{\nu_2})$, whose dimension is $\dim(\GH^{\mathfrak S}_{\nu_1}\cap\GH^{\mathfrak S}_{\nu_2})$, see~\cite[Proposition~3.34]{KM}. Moreover, $\GH^{\mathfrak S}_\nu=J_\nu^\perp$ with respect to the apolarity pairing, so the spaces $\GH^{\mathfrak S}_\nu$ form the dual lattice of the ideals $J_\nu$, and the $(1^N)$-weight space of $S^2V\otimes V^{\otimes(N-2)}$ has dimension $\frac{N!}2$.
	
	\subsection{Set-up}\label{ssec:setup-new}
	
	Throughout the section $N\ge2$ is fixed, $0\le k\le j\le N/2$, and $\lambda=(N-j,j)$. For weights with at most two nonzero entries we write $(a,b):=a\eps_1+b\eps_2$; e.g.\ $\bW_{(N-k,k)}$ is the global Weyl module of highest weight $(N-k)\eps_1+k\eps_2$ with the cyclic vector $w$, and $M[(a,b)]$ denotes the weight space of a module $M$. Let $\bA:=\U(\fh\otimes x\mathbb C[x])$ be the algebra of the Cartan currents of positive degree.
	
	We realize $\mathcal H_j=W_{(N-j)\eps_1}\star_yW_{j\eps_1}$ with the fusion parameters $0$ and $1$: on the tensor product $W_{(N-j)\eps_1}\otimes W_{j\eps_1}$ an element $ax^ry^s$ acts as $\delta_{s,0}(ax^r)\otimes1+1\otimes(ax^r)$, i.e.\ $y$ acts by $0$ on the first and by $1$ on the second factor. Let $v_N:=w_{N-j}\otimes w_j$ be the product of the cyclic vectors, and let $\mathcal F_l$ be the span of the products of elements of $\U(\gl[x,y])$ of total $y$-degree $\le l$ applied to $v_N$. Each $\mathcal F_l$ is a $\gl[x]$-submodule, and
	\[
	\mathcal H_j^k=\mathcal F_k/\mathcal F_{k-1},\qquad v_k:=\text{the class of }(E_{21}y)^kv_N=w_{N-j}\otimes E_{21}^kw_j .
	\]
	Since $\fn_+[x,y]$ kills $v_N$ and $\fh[x,y]$ acts on $v_N$ by scalars, the tensor product equals $\U(\fn_-[x,y])v_N$. A monomial in $\fn_-[x,y]$ of weight $-i\alpha_1$ is a product of $i$ elements $E_{21}x^ay^b$, and $E_{21}x^ay^b$ acts on the tensor product as $E_{21}x^ay$ if $b\ge1$. Hence
	\begin{equation}\label{eq:weight-filtration}
		\big(W_{(N-j)\eps_1}\otimes W_{j\eps_1}\big)[(N-i,i)]\subset\mathcal F_i\qquad(i\ge0).
	\end{equation}
	
	\begin{lemma}\label{lem:y2}
		Let $\rho_j\colon S^N(V[x,y])\twoheadrightarrow\mathcal H_j$, $e_1^N\mapsto v_N$, be the surjection of Theorem~\ref{thm:GH-fusion}, and let $S^N(V[x,y])_k$ be the component of $y$-degree $k$.
		\begin{enumerate}
			\item $\rho_j$ vanishes on the span $I$ of the monomials having a factor in $V[x]y^{\ge2}$.
			\item $S^N(V[x,y])_k/(I\cap S^N(V[x,y])_k)\simeq\bU_k:=S^{N-k}(V[x])\otimes S^k(V[x]y)$, and this $\gl[x]$-module is generated by $u_0:=e_1^{N-k}(e_2y)^k$.
			\item $\rho_j(u_0)$ is a nonzero multiple of $v_k$. In particular, $\mathcal H_j^k=\U(\gl[x])\,v_k$.
		\end{enumerate}
	\end{lemma}
	\begin{proof}
		(1) A monomial with a factor $e_c x^a y^b$, $b\ge2$, equals $E_{cd}x^a y^b$ applied to the monomial in which this factor is replaced by $e_d$, for a colour $d$ which does not occur in it (possible as $m>N$). On the tensor product the operator $E_{cd}x^a y^b$ maps $\mathcal F_l$ to $\mathcal F_{l+1}\subset\mathcal F_{l+b-1}$, and therefore acts by zero from $\mathcal F_l/\mathcal F_{l-1}$ to $\mathcal F_{l+b}/\mathcal F_{l+b-1}$.
		
		(2) The isomorphism is clear. For $m\ge N$ the $\gl[x]$-module $V[x]^{\otimes N}$ is generated by $e_1\otimes\dots\otimes e_N$: the currents $E_{aa}x^r$ multiply the factor of colour $a$ by $x^r$, the elements of weight $0$ of $\U(\gl)$ permute the colours, and every irreducible polynomial $\gl$-module of degree $N$ has a nonzero $(1^N)$-weight space. Hence $\bU_k$ is generated by $e_1e_{c_2}\cdots e_{c_{N-k}}\otimes(e_2y)(e_{d_2}y)\cdots(e_{d_k}y)$ with pairwise distinct colours $1,2,c_i,d_i$, and this vector is obtained from $u_0$ by the lowering operators $E_{c_i1}$ and $E_{d_i2}$.
		
		(3) $(E_{21}y)^ke_1^N=\frac{N!}{(N-k)!}\,u_0$, and $\rho_j$ is a homomorphism. Consequently, thanks to (1) and (2) we see that $\mathcal{H}_j^k$ is generated by $v_k$ as a $\gl[x]$-module.
	\end{proof}

	\subsection{Weight spaces of the local Weyl module $W_{j\eps_1}$}
	
	In this subsection we describe the weight spaces of the local Weyl module $W_{j\eps_1}$ as modules over the Cartan currents; they control the highest weight spaces of the modules $\mathcal H_j^k$.
	
	Let $\nu=\sum_c\nu_c\eps_c$ be a weight of $W_{j\eps_1}$, so that $\sum_c\nu_c=j$. Lemma~\ref{lem:weight-space} below describes the weight space $W_{j\eps_1}[\nu]$ as a module over $\fh[x]$. In order to state it, we introduce appropriate notation:
	\\
	For every colour $c=1,\dots,m$ consider the variables $t^c_1,t^c_2,\dots$, and let $e_r(t^c_\bullet)$, $h_r(t^c_\bullet)$, $p_r(t^c_\bullet)$ be the elementary, complete and power sum symmetric functions in them. Let
	\[
	\psi\colon\Lam(t^1_\bullet)\otimes\dots\otimes\Lam(t^m_\bullet)\to\U(\fh[x]),\qquad \psi\big(p_r(t^c_\bullet)\big):=E_{cc}x^r\quad(r\ge1),
	\]
	be the homomorphism from the tensor product of the rings of symmetric functions in the variables of each colour (a polynomial ring in the $p_r(t^c_\bullet)$) to the commutative algebra of Cartan currents. For the weight $\nu$ we keep only $\nu_c$ variables of colour $c$, i.e.\ we use the restriction $\Lam(t^c_\bullet)\to\mathbb C[t^c_1,\dots,t^c_{\nu_c}]^{\sS_{\nu_c}}$, and put $\sS_\nu=\sS_{\nu_1}\times\dots\times\sS_{\nu_m}$.
	\begin{lemma}\label{lem:weight-space}
		The $\fh[x]$-module $W_{j\eps_1}[\nu]$ is cyclic, generated by $u_\nu:=\prod_{c\ge2}E_{c1}^{\nu_c}w_{j\eps_1}$, and the map $f\mapsto\psi(f)u_\nu$ induces an isomorphism
		\[
		\mathbb C[t^c_l:\ 1\le l\le\nu_c]^{\sS_\nu}\Big/\Big(\sum_{c=1}^mp_r(t^c_\bullet):\ r\ge1\Big)\ \xrightarrow{\ \sim\ }\ W_{j\eps_1}[\nu].
		\]
	\end{lemma}
	\begin{proof}
		Consider the global Weyl module $\bW_{j\eps_1}\simeq S^j(V[x])$. Its weight spaces, as modules over the Cartan currents, are
		\[
		\bW_{j\eps_1}[\nu]\simeq\bigotimes_{c=1}^mS^{\nu_c}\big(e_c\otimes\mathbb C[x]\big)\simeq\mathbb C[t^c_l:\ 1\le l\le\nu_c]^{\sS_\nu},
		\]
		where $E_{cc}x^r$ acts by multiplication by $p_r(t^c_\bullet)$ and the vector $u_\nu$ corresponds to a nonzero constant. The central elements $\sum_{c=1}^mE_{cc}x^r$, $r\ge1$, annihilate the cyclic vector of the local Weyl module and hence the whole module; in particular they annihilate $u_\nu$, and $\psi\big(\sum_cp_r(t^c_\bullet)\big)=\sum_cE_{cc}x^r$. Thus $W_{j\eps_1}[\nu]$ is a quotient of the ring in the lemma. Finally,
		\[
		\dim W_{j\eps_1}[\nu]=\binom{j}{\nu_1,\nu_2,\dots}=\dim\mathbb C[t^c_l]^{\sS_\nu}\Big/\Big(\sum_cp_r(t^c_\bullet):\ r\ge1\Big),
		\]
		the first equality since $W_{j\eps_1}\simeq V^{\otimes j}$ as a $\gl$-module (Theorem~\ref{thm:Weyl-fusion}), and the second one since $\mathbb C[t]^{\sS_\nu}$ is a free module of rank $|\sS_j/\sS_\nu|$ over the symmetric functions in all $j$ variables.
	\end{proof}
	
	\begin{corollary}\label{cor:weight}
		The $\fh[x]$-module $W_{j\eps_1}[(j-k,k)]$ is cyclic, generated by $u:=E_{21}^kw_{j\eps_1}$, and is defined by the relations
		\[
		E_{cc}x^l\,u=0\ (c\geq 3),\qquad (E_{11}+E_{22})x^l\,u=0,\qquad \psi e_r(t^2_{\bullet})\,u=0\ (r>k),\qquad \psi e_s(t^1_{\bullet})\,u=0\ (s>j-k),
		\]
		$l\ge1$. 
		It is isomorphic to $\mathbb C[t^1_1,\dots,t^1_{j-k},t^2_1,\dots,t^2_k]^{\sS_{j-k}\times\sS_k}/(p_r(t^1_\bullet)+p_r(t^2_\bullet):r\ge1)\simeq H^*(\Gr(k,j);\mathbb C)$; in particular its dimension is $\binom jk$ and its graded dimension (with respect to the $x$-degree) is the $q$-binomial coefficient $\qbinom jk$.
	\end{corollary}
	\begin{proof}
		This is Lemma~\ref{lem:weight-space} for $\nu=(j-k,k)$: the kernel of the restriction $\Lam(t^c_\bullet)\to\mathbb C[t^c_1,\dots,t^c_{\nu_c}]^{\sS_{\nu_c}}$ is generated by the $e_s(t^c_\bullet)$ with $s>\nu_c$ (for $c\ge3$, $\nu_c=0$, it is generated by the $p_r(t^c_\bullet)$). The isomorphism with $H^*(\Gr(k,j))$ is Borel's presentation of the cohomology of the Grassmannian~\cite{F}.
	\end{proof}
	
	\begin{lemma}\label{lem:HighestWectorSubquotion2}
		We have the following inclusion for the fusion filtration of the module $W_{(N-j)\varepsilon_1}\otimes W_{j\varepsilon_1}$:
		\[\fn_+[x]\,(E_{21}y)^k v_N \subset \mathcal{F}_{k-1}.\]
		In particular, $w\mapsto v_k$ defines a surjection $\bW_{(N-k,k)}\twoheadrightarrow\mathcal H_j^k$.
	\end{lemma}
	\begin{proof}
		We have $E_{i,i+1}x^a(E_{21}y)^kv_N=0$ for $i\neq 1$,and the vector $E_{12}x^a(E_{21}y)^kv_N$ has weight $(N-k+1,k-1)$, hence it lies in $\mathcal F_{k-1}$ by~\eqref{eq:weight-filtration}.
		Thus the class $v_k$ is annihilated by $\fn_+[x]$. It has weight $(N-k,k)$, the central currents $\sum_cE_{cc}x^r$ act by zero on $\mathcal H_j$ (they do so on both local Weyl modules), and $\mathcal H_j^k$ is a finite-dimensional $\gl$-module. Therefore $v_k$ satisfies the defining relations of $\bW_{(N-k,k)}$, and the map $w\mapsto v_k$ is surjective by Lemma~\ref{lem:y2}(3).
	\end{proof}

	The following lemma describes the highest weight space of $\mathcal H_j^k$.
	\begin{lemma}\label{lem:hw-space}
		The map $u\mapsto$ (the class of $w_{N-j}\otimes u$) is a surjection of $\fh[x]$-modules
		\[W_{j \eps_1}[(j-k,k)] \twoheadrightarrow \mathcal{H}_j^k[(N-k,k)],\qquad E_{21}^kw_j\mapsto v_k.\]
	\end{lemma}
	\begin{proof}
		The map is $\fh[x]$-linear because the Cartan currents of positive degree annihilate $w_{N-j}$, and $w_{N-j}\otimes u\in\mathcal F_k$ by~\eqref{eq:weight-filtration}. The weight space of the tensor product decomposes as
		\[\big(W_{(N-j)\eps_1}\otimes W_{j\eps_1}\big)[(N-k,k)]=\bigoplus_{l=0}^k W_{(N-j)\eps_1}[(N-j-l,l)]\otimes W_{j\eps_1}[(j-k+l,k-l)].\]
		
		By~\eqref{eq:weight-filtration} the weight space of weight $(N-k+1,k-1)$ lies in $\mathcal F_{k-1}$, and so does its image under $E_{21}x^r$, since $\mathcal F_{k-1}$ is a $\gl[x]$-submodule. Moreover, the Weyl modules are cyclic highest weight modules, so
		\[{\rm span}_r\, E_{21}x^r\big(W_{(N-j)\eps_1}[(N-j-l+1,l-1)]\big)=W_{(N-j)\eps_1}[(N-j-l,l)].\]
		Hence, by the Leibniz rule $E_{21}x^ru_1\otimes u_2=E_{21}x^r(u_1\otimes u_2)-u_1\otimes E_{21}x^ru_2$, every summand with $l\ge1$ lies in the summand with $l-1$ modulo $\mathcal F_{k-1}$. By induction on $l$ we obtain
		\[w_{N-j}\otimes W_{j\eps_1}[(j-k,k)]+\mathcal{F}_{k-1}[(N-k,k)]=\mathcal{F}_{k}[(N-k,k)].\]
		This implies the claim.
	\end{proof}
	
	We need the following computation in the universal enveloping algebra.
	\begin{lemma}\label{lem:comm}
		For $c,d\ge0$
		\[
		\psi(e_d(t^1_\bullet))E_{21}x^{c}=\sum_{b=0}^{d}(-1)^{d-b}E_{21}x^{c+d-b}\psi(e_b(t^1_\bullet)).
		\]
	\end{lemma}
	\begin{proof}
		Let $e(z)=\sum_b\psi(e_b(t^1_\bullet))z^b=\exp\big(\sum_{r\ge1}(-1)^{r-1}E_{11}x^rz^r/r\big)$. Since $[E_{11}x^r,E_{21}x^n]=-E_{21}x^{n+r}$, the operator $\mathrm{ad}(E_{11}x^r)$ acts on the span of the $E_{21}x^n$ as $-S^r$, where $S\colon E_{21}x^n\mapsto E_{21}x^{n+1}$. Hence conjugation by $e(z)$ acts there as $\exp\big(-\sum_r(-1)^{r-1}z^rS^r/r\big)=(1+zS)^{-1}$, i.e.
		\begin{equation}\label{eq:conj-e}
			e(z)\,E_{21}x^c\,e(z)^{-1}=\sum_{n\ge0}(-z)^nE_{21}x^{c+n}.
		\end{equation}
		Taking the coefficient of $z^d$ in $e(z)E_{21}x^c=\big(\sum_n(-z)^nE_{21}x^{c+n}\big)e(z)$ we get the claim.
	\end{proof}
	
	\begin{theorem}\label{thm:presentation}
		Let $0\le k\le j\le N/2$. The map $\bW_{(N-k,k)}\to\mathcal H_j^k$, $w\mapsto v_k$, is surjective, and its kernel is the $\U(\gl_V[x])$-submodule generated by the following elements:
		\begin{enumerate}[label=\textup{$(\mathrm{R}\arabic*_j)$}]
			\item \label{R1}  $E_{cc}x^l w$ $(c \geq 3)$, $(E_{11}+E_{22})x^l  w$, $\psi e_r(t^2_{\bullet})w$ $(r>k)$, $\psi e_s(t^1_{\bullet})w$ $(s>j-k)$, where $l\ge1$;
			\item   \label{R2}  $\psi (e_{j-k}(t^1_\bullet))\,E_{21}x^{N-2j}\,a\,w$ for every $a\in\U(\fh[x])$; equivalently (Lemma~\ref{lem:comm})
			\[
			\sum_{b=0}^{j-k}(-1)^{j-k-b}E_{21}x^{N-j-k-b}\psi (e_{b}(t^1_\bullet))\,a\,w .
			\]
		\end{enumerate}
	\end{theorem}
	\begin{proof}
		The surjectivity is Lemma~\ref{lem:HighestWectorSubquotion2}. The relations~\ref{R1} hold in $\mathcal H_j^k$: they hold for $E_{21}^kw_j\in W_{j\eps_1}[(j-k,k)]$ by Corollary~\ref{cor:weight}, and the surjection of Lemma~\ref{lem:hw-space} is $\fh[x]$-linear and sends $E_{21}^kw_j$ to $v_k$. Let us check~\ref{R2}.
		We use the realization of $\mathcal H_j$ from \S\ref{ssec:setup-new}: $\mathcal H_j^k=\mathcal F_k/\mathcal F_{k-1}$ and $v_k$ is the class of $(E_{21}y)^kv_N$. The elements
		$E_{21}x^n$ and $E_{21}x^ny$ act as the coefficients of $z^n$ in the series
		\[
		E_{21}(z)=\sum_n\big(E_{21}x^n\otimes1+1\otimes E_{21}x^n\big)z^n,\qquad
		E_{21}y(z)=\sum_n\big(1\otimes E_{21}x^n\big)z^n ,
		\]
		and $E_{21}(z)-E_{21}y(z)=\sum_n(E_{21}x^n\otimes1)z^n$. All these coefficients commute, and one obviously has
		\begin{equation}\label{eq:filt}
			\text{a product of }
			\left\{\begin{array}{c}
				\text{ at most }l\text{ coefficients of }E_{21}y(z) \text{ and }
				\\
				\text{ an arbitrary amount of coefficients of }E_{21}(z)
			\end{array}
			\right\}
			\text{ applied to }v_N,\text{ lies in }\mathcal F_l .
		\end{equation}
		First, recall, that in a single local Weyl module $W^x_{r\eps_1}$ is a fusion product of $r$ copies of the tautological representations.
		Therefore, for the generating series
		$E_{21}(z)=\sum_{n\geq 0}E_{21}x^nz^n$ we have the following annihilation property:
		\[
		q+m > r \quad \Rightarrow \quad 
		[z^m]\,E_{21}(z)^q\,w_{r\eps_1}=0.
		\]
		Indeed, $W^x_{r\eps_1}$ is the associated graded of $V(\beta_1)\otimes\dots\otimes V(\beta_r)$ (Theorem~\ref{thm:Weyl-fusion}), where $E_{21}(z)=\sum_iE_{21}^{(i)}(1-\beta_iz)^{-1}$ and $(E_{21}^{(i)})^2=0$. Hence $E_{21}(z)^q=q!\sum_{|I|=q}\prod_{i\in I}(1-\beta_iz)^{-1}\prod_{i\in I}E_{21}^{(i)}$, so $\prod_{i=1}^r(1-\beta_iz)\cdot E_{21}(z)^q$ is a polynomial in $z$ of degree $\le r-q$. For $m>r-q$ its $z^m$-coefficient vanishes; this expresses $[z^m]E_{21}(z)^qw$ through $[z^{m-i}]E_{21}(z)^qw$, $i\ge1$, i.e.\ through elements of lower filtration degree, so $[z^m]E_{21}(z)^qw=0$ in the associated graded module.
		Applying this to the two factors of $v_N$ we get
		\begin{align}
			&[z^m]\big(E_{21}(z)-E_{21}y(z)\big)^{k+1}v_N=0\quad (m\ge N-j-k),\tag{a}\\
			&[z^m]\,E_{21}y(z)^qv_N=0\quad (m\ge j-q+1).\tag{b}
		\end{align}
		(For $k=0$ read $E_{21}y(z)^0=1$.)
		
		Recall the generating series $\sum_b\psi e_b(t^1_\bullet)z^b=\exp\big(\sum_{r\ge1}(-1)^{r-1}E_{11}x^rz^r/r\big)$ from the proof of Lemma~\ref{lem:comm},
		and the commutation relations $[E_{11}x^r,E_{21}x^ny^d]=-E_{21}x^{n+r}y^d$ and $E_{11}x^rv_N=0$. Conjugating by this exponential sends
		$E_{21}x^ny^d$ to $\sum_m(-z)^mE_{21}x^{n+m}y^d$ and fixes $v_N$. Hence, 
		\begin{equation}\label{eq:shift}
			\forall n\geq 0 \text{ we have } \psi e_b(t^1_\bullet)\,E_{21}x^n(E_{21}y)^kv_N=(-1)^b\!\!\sum_{m_0+\dots+m_k=b}\!\!E_{21}x^{n+m_0}\,(E_{21}x^{m_1}y)\cdots(E_{21}x^{m_k}y)\,v_N .
		\end{equation}
		Take $b=j-k$ and $n\ge N-2j$. By \eqref{eq:shift} the right-hand side equals
		$$(-1)^{j-k}\sum_{a=n}^{n+j-k}E_{21}x^a\,[z^{n+j-k-a}]\bigl(E_{21}y(z)\bigr)^kv_N.$$
		The terms with $a<n$ may be added: they
		vanish by (b), since then $n+j-k-a>j-k$. Therefore
		\[
		\psi e_{j-k}(t^1_\bullet)\,E_{21}x^n(E_{21}y)^kv_N=(-1)^{j-k}\,[z^{n+j-k}]\big(E_{21}(z)\,E_{21}y(z)^k\big)v_N .
		\]
		Since $n+j-k\ge N-j-k$, relation (a) gives, with $m=n+j-k$,
		\[
		0=[z^m]\big(E_{21}(z)-E_{21}y(z)\big)^{k+1}v_N=\sum_{l=0}^{k+1}(-1)^l\binom{k+1}{l}[z^m]\big(E_{21}(z)^{k+1-l}E_{21}y(z)^l\big)v_N .
		\]
		The term $l=k+1$ vanishes by (b) ($m\ge j-k$). The terms $l\le k-1$ lie in $\mathcal F_{k-1}$ by \eqref{eq:filt}.
		The term $l=k$ is $(-1)^k(k+1)[z^m]\big(E_{21}(z)E_{21}y(z)^k\big)v_N$. Hence,
		\[
		\psi e_{j-k}(t^1_\bullet)\,E_{21}x^n(E_{21}y)^kv_N\in\mathcal F_{k-1}\qquad\text{for all }n\ge N-2j,
		\]
		that is, $\psi e_{j-k}(t^1_\bullet)E_{21}x^{n}v_k=0$ in $\mathcal H_j^k$.
		
		Let us now prove Relation~\ref{R2} for arbitrary $a$. Let $X$ be a polynomial in the coefficients of $E_{21}(z)$ and
		$E_{21}y(z)$. On the vector $Xv_N$ the currents $E_{cc}x^r$ ($c\ge3$) act by $0$, and $E_{11}x^r$ acts as $-E_{22}x^r$. So it suffices to take
		$a$ a product of currents $E_{22}x^r$. We show by induction on the length of $a$ that
		$\psi e_{j-k}(t^1_\bullet)E_{21}x^n\,a\,v_k=0$ in $\mathcal H_j^k$ for all $n\ge N-2j$. The case $a=1$ is what we already explained. For
		$a=E_{22}x^r\,a'$ use $[E_{22}x^r,E_{21}x^n]=E_{21}x^{n+r}$:
		\[
		\psi e_{j-k}(t^1_\bullet)E_{21}x^nE_{22}x^ra'v_k=E_{22}x^r\,\psi e_{j-k}(t^1_\bullet)E_{21}x^na'v_k-\psi e_{j-k}(t^1_\bullet)E_{21}x^{n+r}a'v_k ,
		\]
		and both terms vanish in $\mathcal H_j^k$ by the induction hypothesis (the first because $\mathcal F_{k-1}$ is $\gl[x]$-stable).
		For $n=N-2j$ this is~\ref{R2}; its second form follows from Lemma~\ref{lem:comm}.

		Let us show that relations $(R1_j, R2_j)$ are complete. Let $Q$ be the quotient of $\bW_{(N-k,k)}$ by the submodule generated by~\ref{R1} and~\ref{R2}; we have shown that $Q\twoheadrightarrow\mathcal H_j^k$. By~\ref{R1} and Corollary~\ref{cor:weight}, $\bA w\subset Q$ is a quotient of $W_{j\eps_1}[(j-k,k)]$, so there are homogeneous $a_1,\dots,a_r\in\bA$, $r\le\binom jk$, of degrees $d_1,\dots,d_r$, such that $a_1w,\dots,a_rw$ is a basis of $\bA w$. Since $\fn_+[x]\bA w=0$, the PBW theorem gives $Q=\sum_i\U(\fn_-[x])a_iw$. Consider the filtration of $Q$ by the submodules $Q_{\ge d}:=\sum_{d_i\ge d}\U(\gl[x])a_iw$ and the subquotients $M_i:=\U(\gl[x])a_iw/(\U(\gl[x])a_iw\cap Q_{\ge d_i+1})$. The image of $a_iw$ in $M_i$ is annihilated by $\fn_+[x]$ and by $\fh\otimes x\mathbb C[x]$ (because $(\fh\otimes x^t)\,a_iw\in\bA_{>d_i}w\subset Q_{\ge d_i+1}$ for $t\ge1$), has weight $(N-k,k)$, and generates an integrable module; so $M_i$ is a quotient of the local Weyl module $W_{(N-k,k)}$. Moreover, in the second form of~\ref{R2} for $a=a_i$ all the terms with $b\ge1$ lie in $Q_{\ge d_i+1}$, so $E_{21}x^{N-j-k}a_iw=0$ in $M_i$. By Corollary~\ref{cor:two-rows-presentation}, $M_i$ is a quotient of $W_{(N-k,k)}/\U(\gl[x])E_{21}x^{N-j-k}w\simeq\GPD_{j,k}$, and therefore
		\[
		\dim\mathcal H_j^k\le\dim Q\le\sum_{i=1}^r\dim M_i\le\binom jk\dim\GPD_{j,k}.
		\]
		Summing over $k$ and using Example~\ref{ex:GPD-two-rows} we get 
		$$m^N=\dim\mathcal H_j=\sum_k\dim\mathcal H_j^k\le\sum_k\binom jk\dim\GPD_{j,k}=m^N.$$
		Hence all the inequalities are equalities; in particular $Q\simeq\mathcal H_j^k$, i.e.\ \ref{R1} and~\ref{R2} are the defining relations of $\mathcal H_j^k$. Moreover $r=\binom jk$ and $M_i\simeq\GPD_{j,k}$ for all $i$.
	\end{proof}

	\begin{remark}
		Modulo the Cartan currents of positive degree the relation~\ref{R2} becomes $E_{21}x^{N-j-k}w=0$, the defining relation of $\GPD_{j,k}$ (Corollary~\ref{cor:two-rows-presentation}). Thus~\ref{R2} is a deformation of this relation by the Cartan currents. In particular,
		\begin{itemize}[itemsep=0pt,topsep=0pt]
			\item for $k=j$ the relations~\ref{R1} give $\bA w=0$, and $\mathcal H_j^j\simeq W_{(N-j,j)}=\GPD_{j,j}$; 
			\item
			for $k=0$ we get $\mathcal H_j^0\simeq W_{(N,0)}/\U(\gl[x])E_{21}x^{N-j}w=\GPD_{j,0}$.
		\end{itemize}
	\end{remark}
	
	\begin{corollary}\label{cor:A}
		$\mathcal H_j^k$ has a filtration by $\gl_V[x]$-submodules whose $\binom jk$ subquotients are isomorphic (up to shifts of the $x$-grading) to $\GPD_{j,k}$; thus 
		$$\ch_q\mathcal H_j^k=\qbinom jk\,\ch_q\GPD_{j,k}.$$
	\end{corollary}
	\begin{proof}
		This was shown at the end of the proof of Theorem~\ref{thm:presentation}: the subquotients of the filtration $Q_{\ge d}$ are $M_i\simeq\GPD_{j,k}$, shifted by $d_i$, and $\sum_iq^{d_i}=\qbinom jk$ is the graded dimension of $\bA w\simeq W_{j\eps_1}[(j-k,k)]$ (Corollary~\ref{cor:weight}).
	\end{proof}
	
	\begin{corollary}\label{thm:B}
		The surjection of Lemma~\ref{lem:hw-space} is an isomorphism. Consequently the algebra $\U(\fh[x])/\Ann(v_k)$ of the Cartan currents acting on the highest weight vector $v_k\in\mathcal H_j^k$ is isomorphic to $H^*(\Gr(k,j);\mathbb C)$ (Corollary~\ref{cor:weight}).
	\end{corollary}
	\begin{proof}
		Since $\mathcal H_j^k=\U(\gl[x])v_k=\U(\fn_-[x])\bA v_k$, its highest weight space is $\bA v_k$, of dimension $r=\binom jk=\dim W_{j\eps_1}[(j-k,k)]$ by the proof of Theorem~\ref{thm:presentation}.
	\end{proof}
	
	\begin{remark}\label{thm:Schur-lattice}
		Since $E_{11}x^r$ acts on $v_k$ as $-E_{22}x^r$, the Cartan part of the relations~\ref{R1} is the ideal $(e_r:\ r>k)+(h_s:\ s>j-k)$ of $\Lam=\Lam(t^2_\bullet)$; here we used $e_s(t^1_\bullet)\mapsto(-1)^sh_s(t^2_\bullet)$ under $p_r(t^1_\bullet)\mapsto-p_r(t^2_\bullet)$. This ideal equals $\Span\{s_\mu:\ \mu\not\subset((j-k)^k)\}$, and the quotient by it is $H^*(\Gr(k,j))$~\cite{F}. Thus, for fixed $k$, the Cartan relations are spanned by subsets of the basis of Schur functions and form a decreasing chain in $j$; this is the chain $(G_j)$ of Proposition~\ref{thm:two-chains}.
	\end{remark}
	
	For a fixed $k$ we now compare the relations for different $j$. For $j<k$ we have $\mathcal H_j^k=0$; accordingly, the relations~\ref{R1} contain $\psi e_0(t^1_\bullet)w=w$ in this case.
	
	\begin{lemma}\label{lem:ChainsOfRelations}\label{lem:lift}
		Let $G_j,F_j\subset\bW_{(N-k,k)}$ be the submodules generated by the relations~\ref{R1} and~\ref{R2} respectively (so that $G_j+F_j$ is the kernel of $\bW_{(N-k,k)}\twoheadrightarrow\mathcal H_j^k$). Then for $j<\lfloor N/2\rfloor$
		\[
		G_{j+1}\subset G_j,\qquad F_j\subset G_j+F_{j+1}.
		\]
		In words: the relations $(\mathrm{R1}_j)$ imply the relations $(\mathrm{R1}_{j+1})$, and the relations $(\mathrm{R2}_{j+1})$ together with $(\mathrm{R1}_j)$ imply $(\mathrm{R2}_j)$.
	\end{lemma}
	\begin{proof}
		The relations $(\mathrm{R1}_j)$ differ from $(\mathrm{R1}_{j+1})$ by one extra relation $\psi e_{j-k+1}(t^1_{\bullet})w=0$.
		Denote $e_d:=\psi e_d(t^1_\bullet)$. 
		Comparing the coefficients of $z^{d+1}$ in $e(z)E_{21}x^{n}=E_{21}x^{n}e(z)-z\,e(z)E_{21}x^{n+1}$, which is~\eqref{eq:conj-e}, we get 
		$$e_d E_{21}x^{n+1}=E_{21}x^{n} e_{d+1} -  e_{d+1} E_{21}x^{n} \text{ and }
		e_{d+1}E_{21}x^{n}=E_{22}x\,e_{d+1}E_{21}x^{n-1}-e_{d+1}E_{21}x^{n-1}E_{22}x,$$ because $[E_{22}x,E_{21}x^{n-1}]=E_{21}x^n$ and $E_{22}x$ commutes with $e_{d+1}$. For $d=j-k$ and $n=N-2j-1$ we obtain, for every $a\in\U(\fh[x])$,
		\begin{multline}\label{eq:lift}
			e_{j-k}E_{21}x^{N-2j}a\,w = \\ E_{21}x^{N-2j-1}a\,e_{j-k+1}w-E_{22}x\,\big(e_{j-k+1}E_{21}x^{N-2(j+1)}a\,w\big)+e_{j-k+1}E_{21}x^{N-2(j+1)}\big(E_{22}x\,a\big)w .
		\end{multline}
		The first term lies in $G_j$, and the last two terms lie in $F_{j+1}$.
	\end{proof}

	\begin{lemma}\label{lem:GPD-chain}
		In $W=W_{(N-k,k)}$ put $R_s=\U(\gl[x])\,E_{21}x^sw$. Then $R_{s+1}\subset R_s$. Hence for $k\le j<\lfloor N/2\rfloor$ one has surjections $\GPD_{j,k}\twoheadrightarrow\GPD_{j+1,k}$, the kernels $R_{N-j-k}$ form a chain, and the lattice of submodules of $W$ generated by them is a chain, in particular distributive. In the global module the same inclusion holds modulo the Cartan relations \cbl{(Lemma~\ref{lem:ChainsOfRelations})}.
	\end{lemma}
	\begin{proof}
		In $W$ we have $E_{22}x\,w=0$, hence $E_{21}x^{s+1}w=[E_{22}x,E_{21}x^s]w=E_{22}x\cdot E_{21}x^sw\in R_s$. The surjections follow from Corollary~\ref{cor:two-rows-presentation}, since $N-j-k$ decreases when $j$ grows.
	\end{proof}

	\subsection{The main theorems in the $\gl$-form and their proofs}\label{ssec:proofs}
	
	In this subsection we provide the proofs of the Schur--Weyl dual forms of Theorems~\ref{thm:lattice} and~\ref{thm:N!/2}; for the equivalence see the discussion after Theorem~\ref{thm:N!/2} and~\cite[\S3.3]{KM}. As before, $\cJ_{(N-j,j)}=\ker\rho_j\subset S^N(V[x,y])$.

	First, let us recall the following classical theorem on distributive lattices:
	\begin{lemma}[Birkhoff]\label{lem:birkhoff}
		In the lattice of subspaces of a finite-dimensional vector space, two chains $A_0\supseteq\dots\supseteq A_p$ and $B_0\subseteq\dots\subseteq B_p$ generate a distributive sublattice. In particular, if $K_j=A_j+B_j$ for all $j$, then the $K_j$ generate a distributive lattice.
	\end{lemma}
	\begin{proof}
		Two flags admit a common adapted basis. Every $A_j$ and $B_j$ is spanned by a subset of this basis, and such subspaces form a Boolean lattice. This is Birkhoff's theorem on modular lattices generated by two chains~\cite{Birkhoff}.
	\end{proof}

	Fix $k$ and let $K_j:=G_j+F_j=\ker\big(\bW_{(N-k,k)}\twoheadrightarrow\mathcal H_j^k\big)$ for $0\le j\le\lfloor N/2\rfloor$ (In particular, $K_j=G_j=\bW_{(N-k,k)}$ for $j<k$).
	
	\begin{proposition}\label{thm:two-chains}
		Put $B_j:=\bigcap_{j'\ge j}K_{j'}$, so that $B_j\subset B_{j+1}$. Then $K_j=G_j+B_j$ for all $j$. Thus every $K_j$ is the sum of a member of the decreasing chain $(G_j)$ and a member of the increasing chain $(B_j)$.
	\end{proposition}
	\begin{proof}
		Clearly $G_j+B_j\subset K_j$. We prove the equality by downward induction on $j$; for $j=\lfloor N/2\rfloor$ we have $B_j=K_j$. Let $j<\lfloor N/2\rfloor$ and assume $K_{j+1}=G_{j+1}+B_{j+1}$. By Lemma~\ref{lem:ChainsOfRelations},
		\[
		K_j=G_j+F_j\subset G_j+F_{j+1}\subset G_j+K_{j+1}=G_j+G_{j+1}+B_{j+1}=G_j+B_{j+1}.
		\]
		Since $G_j\subset K_j$, the modular law gives $K_j=K_j\cap(G_j+B_{j+1})=G_j+(K_j\cap B_{j+1})=G_j+B_j$.
	\end{proof}
	
	\begin{theorem}\label{thm:lattice-gl}
		Let $N\ge2$. The $\gl_V[x,y]$-submodules $\cJ_{(N-j,j)}\subset S^N(V[x,y])$, $0\le j\le\lfloor N/2\rfloor$, generate a distributive lattice.
	\end{theorem}	
	\begin{proof}
		Let $\cJ_j:=\cJ_{(N-j,j)}=\ker\rho_j$. The submodules $\cJ_j\subset S^N(V[x,y])$ are bigraded, so it suffices to show that for every $k$ the components $\cJ_j\cap S^N(V[x,y])_k$ generate a distributive lattice. By Lemma~\ref{lem:y2} all of them contain $I\cap S^N(V[x,y])_k$, so we may pass to the quotient $\bU_k=\U(\gl[x])u_0$, in which they become the kernels of $\bU_k\twoheadrightarrow\mathcal H_j^k$, $u_0\mapsto c\,v_k$ ($c\ne0$).
		
		Thanks to Birkhoff Lemma~\ref{lem:birkhoff} and Proposition~\ref{thm:two-chains} we see, that the desired subspaces generate a distributive lattice of submodules in $\bW_{(N-k,k)}$. So it remains to  pass from the modules $\bW_{(N-k,k)}$ to $S^N(V[x,y])$ and we use the following elementary observation for that (applied to $M=\bU_k$, $M'=\bW_{(N-k,k)}$, $m=c^{-1}u_0$, $m'=w$).
		
		Let $M=\U m$ be a cyclic module over an algebra $\U$. The map $L\mapsto Lm$ is an isomorphism from the lattice of left ideals $L\supseteq\Ann(m)$ to the lattice of submodules of $M$. Hence, if $M=\U m$ and $M'=\U m'$ are cyclic and $K_i\subset M$, $K'_i\subset M'$ are submodules with $M/K_i\simeq M'/K'_i$ (isomorphisms sending $m$ and $m'$ to the same vector) for all $i$, then $K_i\mapsto K'_i$ extends to an isomorphism of the lattices generated by the $K_i$ and by the $K'_i$: both are isomorphic to the lattice generated by the left ideals $\Ann(m\bmod K_i)=\Ann(m'\bmod K'_i)$.
		
		Consequently, the lattice generated by the submodules $\cJ_j\cap S^N(V[x,y])_k$ is isomorphic to the lattice generated by $K_0,\dots,K_{\lfloor N/2\rfloor}\subset\bW_{(N-k,k)}$. All these submodules are graded by the weights and the $x$-degree, with finite-dimensional components. In each component, Proposition~\ref{thm:two-chains} writes $K_j=G_j+B_j$ with a decreasing chain $(G_j)$ and an increasing chain $(B_j)$, so the $K_j$ generate a distributive lattice by Lemma~\ref{lem:birkhoff}. Hence the $\cJ_j$ generate a distributive lattice of $\gl_V[x,y]$-submodules.
	\end{proof}
	
	\begin{corollary}\label{cor:common-basis}
		Let $N\ge2$. There is a basis of $\mathbb C[\bx,\by]$ consisting of bihomogeneous polynomials such that every Garsia--Haiman module $\GH^{\mathfrak S}_\lambda$ with $\lambda\vdash N$, $\ell(\lambda)\le2$, is spanned by a subset of this basis. Similarly, there is a basis of $S^N(V[x,y])$, homogeneous with respect to the weights and the bidegree, such that every submodule $\cJ_{(N-j,j)}$, $0\le j\le\lfloor N/2\rfloor$, is spanned by a subset of it.
	\end{corollary}
	\begin{proof}
		All the spaces involved are bigraded, with finite-dimensional bigraded components (for $S^N(V[x,y])$ we also split into weight spaces), so it suffices to construct such a basis in every component. By Theorem~\ref{thm:lattice} (resp.\ Theorem~\ref{thm:lattice-gl}) the subspaces $\GH^{\mathfrak S}_{(N-j,j)}$ (resp.\ $\cJ_{(N-j,j)}$) generate a distributive lattice. A finite collection of subspaces of a finite-dimensional vector space generates a distributive lattice if and only if there is a basis such that each of the subspaces is spanned by a subset of it~\cite[Chapter~1, \S7]{PP} (for two chains this is Lemma~\ref{lem:birkhoff}).
	\end{proof}	
	An explicit basis of the Garsia--Haiman modules of two-column diagrams (equivalently, after the exchange $x\leftrightarrow y$, of two-row diagrams), the ``kicking basis'', was constructed by Assaf and Garsia~\cite{Assaf-Garsia-two_rows}, by means of orbit harmonics. It is constructed separately for each diagram, and~\cite{Assaf-Garsia-two_rows} does not discuss its compatibility with the other diagrams of the same size. We do not know an explicit common basis as in Corollary~\ref{cor:common-basis}: the corollary only asserts its existence.		
	\begin{theorem}\label{thm:N!/2-gl}
		Let $1\le j\le N/2$, $\nu_1=(N-j,j)$ and $\nu_2=(N-j+1,j-1)$. Then the common quotient $S^N(V[x,y])/(\cJ_{\nu_1}+\cJ_{\nu_2})$ is isomorphic to $S^2V\otimes V^{\otimes(N-2)}$ as a $\gl_V$-module; equivalently, its character is $h_2p_1^{N-2}$.
	\end{theorem}
	\begin{proof}
		Let $C:=S^N(V[x,y])/(\cJ_j+\cJ_{j-1})$ be the common quotient, where $\cJ_j=\cJ_{(N-j,j)}$.
		
		By~\cite[Theorem~B]{KM}, $C$ surjects onto a module isomorphic to $S^2V\otimes V^{\otimes(N-2)}$ as a $\gl_V$-module and this gives the lover bound of $N!/2$.
		
		To support the upper bound of $N!/2$ we follow the proof of Theorem~\ref{thm:lattice-gl}. Indeed, the $y^k$-component of $C$ is isomorphic to $\bW_{(N-k,k)}/(K_j+K_{j-1})$, which is a quotient of $\bW_{(N-k,k)}/(G_{j-1}+F_j)$, i.e.\ of the module defined by the relations $(\mathrm{R1}_{j-1})$ and $(\mathrm{R2}_j)$. For $k\ge j$ it vanishes, since then $G_{j-1}=\bW_{(N-k,k)}$. For $k\le j-1$ the argument at the end of the proof of Theorem~\ref{thm:presentation} applies verbatim: by $(\mathrm{R1}_{j-1})$ the space $\bA w$ is a quotient of $W_{(j-1)\eps_1}[(j-1-k,k)]$, of dimension $\binom{j-1}k$, and by $(\mathrm{R2}_j)$ the corresponding subquotients are quotients of $\GPD_{j,k}$. Hence, coefficientwise,
		\[
		\ch C\le\sum_{k=0}^{j-1}\binom{j-1}k\ch\GPD_{j,k}=h_2\,p_1^{N-2j}\sum_{k=0}^{j-1}\binom{j-1}kh_2^{\,j-1-k}e_2^k=h_2\,p_1^{N-2j}(h_2+e_2)^{j-1}=h_2\,p_1^{N-2}.
		\]
		Since $h_2p_1^{N-2}=\ch\big(S^2V\otimes V^{\otimes(N-2)}\big)$, the surjection of the lower bound is an isomorphism of $\gl_V$-modules.
	\end{proof}
	
	\appendix
	
	\section{Symmetric functions and plethysm}\label{app:symm}
	
	\subsection{Definitions}\label{ssec:symm-def}
	In this appendix we recall the properties of symmetric functions that are used in the paper. Denote
	\[
	\Lam_k:=\mathbb{C}[x_1,\dots,x_k]^{\mathfrak S_k},\qquad \Lam_k(q):=\mathbb{C}(q)[x_1,\dots,x_k]^{\mathfrak S_k}=\Lam_k\otimes_{\mathbb{C}}\mathbb{C}(q).
	\]
	The ring $\Lam_k$ is graded by the degree, $\Lam_k=\bigoplus_{n\ge0}\Lam_k^n$. For $m>k$ there is the map $\rho_{m,k}\colon\Lam_m\to\Lam_k$, $x_{k+1},\dots,x_m\mapsto0$; clearly $\rho_{m,k}(\Lam_m^n)=\Lam_k^n$. The graded ring of symmetric functions in infinitely many variables is the inverse limit of the rings $\Lam_k^n$ in each degree:
	\[
	\Lam^n:=\varprojlim_k\Lam_k^n,\qquad\Lam:=\bigoplus_{n\ge0}\Lam^n,\qquad \Lam(q):=\Lam\otimes_\mathbb{C}\mathbb{C}(q).
	\]
	There are natural surjective homomorphisms $\rho_k\colon\Lam\to\Lam_k$, which are isomorphisms in degrees $\le k$.
	
	The rings $\Lam_k$ and $\Lam$ are polynomial rings with several standard sets of generators. One is the set of Newton sums $p_l(X)=\sum_ix_i^l$:
	\[
	\Lam=\mathbb{C}[p_1,p_2,\dots],\qquad \Lam_k=\mathbb{C}[\rho_k(p_1),\dots,\rho_k(p_k)].
	\]
	Another is the set of the complete symmetric functions $h_l(X)=\sum_{a_1+a_2+\dots=l}\prod_ix_i^{a_i}$ (and, similarly, of the elementary symmetric functions $e_l$; recall $\sum_lh_lz^l=\exp\sum_l\frac{p_lz^l}{l}$ and $\sum_le_lz^l=\exp\sum_l\frac{(-1)^{l-1}p_lz^l}l$).
	
	A basis of $\Lam^n$ is given by the Schur functions. For a partition $\kappa=(\kappa_1\ge\dots\ge\kappa_k)\vdash n$ put
	\[
	s_\kappa(x_1,\dots,x_k):=\det\big(x_i^{\kappa_l+k-l}\big)_{i,l=1}^k\Big/\det\big(x_i^{k-l}\big)_{i,l=1}^k .
	\]
	The Littlewood--Richardson coefficients $c^\kappa_{\kappa',\kappa''}$ are defined by $s_{\kappa'}s_{\kappa''}=\sum_\kappa c^{\kappa}_{\kappa'\kappa''}s_\kappa$.
	
	\begin{proposition}\label{prop:Schur-basis}
		The functions $s_\kappa(x_1,\dots,x_k)$, $|\kappa|=n$, $\ell(\kappa)\le k$, form a basis of $\Lam_k^n$. Moreover $\rho_{k+1,k}s_\kappa(x_1,\dots,x_{k+1})=s_\kappa(x_1,\dots,x_k)$ if $\ell(\kappa)\le k$ and $0$ otherwise.
	\end{proposition}
	Therefore $s_\kappa(X):=\varprojlim s_\kappa(x_1,\dots,x_k)$ is well defined, and the $s_\kappa$, $|\kappa|=n$, form a basis of $\Lam^n$; for $\kappa=(n)$ we get $h_n$.
	
	\subsection{Plethysm}\label{ssec:plethysm}
	Let $\pi$ be the endomorphism of the ring $\Lam((q))$ of formal Laurent series (rational functions of $q$ are identified with series in $q$) given on the generators by
	\[
	\pi(p_l(X)):=\frac{p_l(X)}{1-q^l}.
	\]
	For $f\in\Lam$ we denote $f\big[\tfrac{X}{1-q}\big]:=\pi(f)$, and for the inverse map $f\big[(1-q)X\big]:=\pi^{-1}(f)$, $\pi^{-1}(p_l)=(1-q^l)p_l$. This is the plethystic substitution $X\mapsto X/(1-q)$: it replaces every variable $x_i$ by the formal infinite sum $x_i(1+q+q^2+\dots)$, i.e.\ $p_l$ by $p_l(x_i q^{a})_{i,a}$.
	
	\begin{lemma}\label{lem:qlog}
		\[
		\exp\Big(\sum_{i\ge1}\frac{Y^i}{i(1-q^i)}\Big)=\sum_{i\ge0}\frac{Y^i}{\prod_{j=1}^i(1-q^j)} .
		\]
	\end{lemma}
	\begin{proof}
		$\sum_{i\ge1}\frac{Y^i}{i(1-q^i)}=\sum_{a\ge0}\sum_{i\ge1}\frac{(q^aY)^i}i=\sum_{a\ge0}\ln\frac1{1-q^aY}$, hence the left-hand side is $\prod_{a\ge0}\frac1{1-q^aY}$, which equals the right-hand side by the $q$-binomial theorem (Euler's identity)~\cite[Theorem~2.1]{A}.
	\end{proof}
	
	\begin{lemma}\label{lem:pleth-h}
		\[
		h_n\Big[\frac X{1-q}\Big]=\frac1{(1-q)\cdots(1-q^n)}\sum_{\substack{a_1,a_2,\dots\in\mathbb N\\ \sum a_i=n}}\binom n{a_1,a_2,\dots}_q\prod_ix_i^{a_i},
		\]
		where $\binom n{a_1,a_2,\dots}_q$ is the $q$-multinomial coefficient.
	\end{lemma}
	\begin{proof}
		Let $z$ be a formal variable. Then $\sum_{n}h_n[\tfrac X{1-q}]z^n=\exp\sum_{i}\frac{p_i(X)}{i(1-q^i)}z^i=\prod_k\exp\sum_i\frac{(x_k z)^i}{i(1-q^i)}$, and by Lemma~\ref{lem:qlog} this is $\prod_k\sum_{a\ge0}\frac{x_k^a z^a}{(q)_a}$. Collecting the coefficients of $z^n$ and using $\prod_i\frac1{(q)_{a_i}}=\frac1{(q)_n}\binom n{a_1,a_2,\dots}_q$, $(q)_n:=\prod_{j=1}^n(1-q^j)$, we obtain the claim.
	\end{proof}
	
	\section{Categories of representations of current algebras and plethysm}\label{app:categories}
	
	The main goal of this appendix is to explain the nature of the plethysm in the representation theory of current algebras: the plethystic substitution $X\mapsto X/(1-q)$ is the character of the functor $V\mapsto V[x]=V\otimes\mathbb{C}[x]$, and the Schur functors of $V[x]$ are the projective modules in the block of polynomial representations of $\gl_m[x]$. This is the representation-theoretic ($\gl_V$-) side of the Schur--Weyl duality for which the plethysm of Haiman's theory appears on the symmetric group side.
	
	\subsection{Weyl modules}\label{ssec:app-weyl}
	Let $\fg$ be a simple Lie algebra (we will also consider $\fg=\gl_m$) and $\fg[x]=\fg\otimes\mathbb{C}[x]$. For a dominant integral weight $\lambda$ let $W_\lambda$ be the local and $\bW_\lambda$ the global Weyl module (see~\cite{CFK}): $\bW_\lambda$ is the cyclic $\fg[x]$-module with the cyclic vector $w_\lambda$ of $\fh$-weight $\lambda$ and the relations
	\[
	(\fn_+\otimes\mathbb{C}[x])\,w_\lambda=0,\qquad(e_{-\alpha}\otimes1)^{\langle\lambda,\alpha^\vee\rangle+1}w_\lambda=0\ \ (\alpha>0),
	\]
	and $W_\lambda$ is its quotient by the additional relation $(\fh\otimes x\mathbb{C}[x])w_\lambda=0$. Both modules are graded by the weights and by the $q$-degree $\deg_qw_\lambda=0$, $\deg_q(a\otimes x^i)=i$ for $a\in\fg$; each graded component is a $\fg$-module. For a graded $\fg$-module $M=\bigoplus_kM[k]$ with finite dimensional components we put $\ch_qM=\sum_kq^k\ch M[k]$. 
	
	Let $\bA=\U(\fh\otimes x\mathbb{C}[x])$. For $a\in\bA$ the submodule $\U(\fg[x])a\,w_\lambda\subset\bW_\lambda$ is cyclic and satisfies the defining relations of $\bW_\lambda$; therefore $\bA$ acts on $\bW_\lambda$ by endomorphisms ($a\cdot uw_\lambda:=u\,aw_\lambda$; since $\bA$ is commutative this is an action of $\bA$ rather than of $\bA^{op}$) commuting with $\fg[x]$. Let $\bA_\lambda=\bA/\ker\varphi$, $\varphi(a)=aw_\lambda$; $\bA_\lambda$ is isomorphic to the space of vectors of $\fh$-weight $\lambda$ in $\bW_\lambda$. We have the following properties (see~\cite{CFK,FL2,N}):
	\begin{itemize}
		\item $W_\lambda$ is finite dimensional, $\dim W_\lambda=\prod_{i}(\dim W_{\omega_i})^{\langle\lambda,\alpha_i^\vee\rangle}$, and $\U(\fg)w_\lambda\simeq V_\lambda\subset W_\lambda$;
		\item the action of $\bA_\lambda$ on $\bW_\lambda$ is free, and if $\lambda=\sum_im_i\omega_i$, then $\bA_\lambda\simeq\bigotimes_i\mathbb{C}[x_{i,1},\dots,x_{i,m_i}]^{\mathfrak S_{m_i}}$ with $h_{\alpha_i}x^k\mapsto x_{i,1}^k+\dots+x_{i,m_i}^k$;
		\item $\ch_q\bW_\lambda=\ch_qW_\lambda\prod_i(q)^{-1}_{\langle\lambda,\alpha_i^\vee\rangle}$, and $W_\lambda\simeq\bW_\lambda/\bA_\lambda^+\bW_\lambda$ ($\bA_\lambda^+$ the augmentation ideal);
		\item $\bW_{N\omega_1}\simeq S^N(V\otimes\mathbb{C}[x])$ for $\mathfrak{sl}_m[x]$.
	\end{itemize}
	Finally, $P_\lambda$ denotes the projective cover of the irreducible module $V_\lambda$ in the category of graded $\fg[x]$-modules with finite dimensional graded components (see~\cite{BBCKL,CI}); explicitly $P_\lambda=\mathrm{Ind}_\fg^{\fg[x]}V_\lambda\simeq\U(\fg\otimes x\mathbb{C}[x])\otimes V_\lambda$ (this is for the full category of graded modules; in the block considered below the projective cover is smaller).
	
	\subsection{The $\gl_m$ Weyl modules}\label{ssec:app-gl}
	Let $\lambda=(\lambda_1,\dots,\lambda_m)$ be a partition, regarded as a dominant weight $\sum\lambda_i\eps_i$ of $\gl_m$, and $\bar\lambda=\sum_{i=1}^{m-1}(\lambda_i-\lambda_{i+1})\omega_i$ the corresponding $\mathfrak{sl}_m$-weight. Since $\gl_m[x]=\mathfrak{sl}_m[x]\oplus\mathrm{Id}[x]$ with the central $\mathrm{Id}[x]=\mathrm{Id}\otimes\mathbb{C}[x]$, we define
	\[
	\bW_\lambda:=\bW_{\bar\lambda}\otimes\U(\mathrm{Id}[x])/\ker\varphi_{\lambda_m},\qquad\varphi_k\colon\mathrm{Id}[x]\to\mathbb{C}[x_1,\dots,x_k],\ \ \mathrm{Id}\otimes x^i\mapsto x_1^i+\dots+x_k^i,
	\]
	where $\mathfrak{sl}_m[x]$ acts on the first factor and $\mathrm{Id}[x]$ on the second. For $\lambda_m=0$ (polynomial weights) the elements $\mathrm{Id}\otimes x^i$, $i\ge1$, act by zero; this is the module used in Section~\ref{sec:fusion}.
	
	\subsection{Categories of polynomial representations}\label{ssec:app-categories}
	
	\begin{definition}
		$\Rep(\gl_m[x])^{N\Box}$ is the category of integrable $\gl_m[x]$-modules all of whose $\gl_m$-highest weights are partitions of $N$.
	\end{definition}
	
	\begin{lemma}
		$\Rep(\gl_m[x])^{N\Box}$ is a block of the category of integrable modules, and $\Rep^{\mathrm{poly}}(\gl_m[x])=\bigoplus_{N\ge0}\Rep(\gl_m[x])^{N\Box}$.
	\end{lemma}
	\begin{proof}
		The central element $h=\mathrm{Id}\in\gl_V$ gives $h\otimes1$ in the center of $\gl_m[x]$, and $hv=Nv$ for every $v\in V_\lambda$, $\lambda\vdash N$. So $\Rep(\gl_m[x])^{N\Box}$ is the eigensubcategory of $h\otimes1$ with the eigenvalue $N$.
	\end{proof}
	
	\begin{lemma}
		If $m>N$, then $\Rep(\gl_m[x])^{N\Box}$ is equivalent to the category $\Rep(\mathfrak{sl}_m[x])^{\le N\omega_1}$ of modules whose highest weights are $\le N\omega_1$.
	\end{lemma}
	\begin{proof}
		The restriction $\mathrm{Res}\colon\Rep(\gl_m[x])\to\Rep(\mathfrak{sl}_m[x])$ maps $\Rep(\gl_m[x])^{N\Box}$ to $\Rep(\mathfrak{sl}_m[x])^{\le N\omega_1}$. The inverse functor equips an $\mathfrak{sl}_m[x]$-module with the action of $\mathrm{Id}\otimes\mathbb{C}[x]$ by $(\mathrm{Id}\otimes1)v=Nv$, $(\mathrm{Id}\otimes x^r)v=0$, $r\ge1$.
	\end{proof}
	
	For $N_1+N_2=N$, $\lambda'\vdash N_1$, $\lambda''\vdash N_2$ one has $V_{\lambda'}\otimes V_{\lambda''}\simeq\bigoplus_{\lambda\vdash N}V_\lambda^{\oplus c^\lambda_{\lambda'\lambda''}}$, hence $$\Rep(\gl_m[x])^{N_1\Box}\otimes\Rep(\gl_m[x])^{N_2\Box}\subset\Rep(\gl_m[x])^{N\Box}.$$ 
	The modules of $\Rep(\gl_m[x])^{N\Box}$ are graded by the weights of $\gl_m$ and the degree of $x$, and the character is a series $\ch M=\sum_iq^i\ch(M)_i$ with $\ch(M)_i\in\mathbb{C}[X_1,\dots,X_m]^{\mathfrak S_m}$ homogeneous of degree $N$; $\ch(M_1\otimes M_2)=\ch M_1\ch M_2$.
	
	\subsection{Projective modules and plethysm}\label{ssec:app-projective}
	
	By $V[x]=V\otimes\mathbb{C}[x]$ we denote the free $\gl_m[x]$-module of rank one on the vector representation; its graded character is $\ch V[x]=(X_1+\dots+X_m)/(1-q)$. For a graded vector space $U$ and a partition $\kappa$ the character of the Schur functor $S_\kappa(U)$ is the plethysm $s_\kappa[\ch U]$ (this is the definition of the plethysm; in particular $p_l[\ch U]$ is obtained from $\ch U$ by the substitution $X_i\mapsto X_i^l$, $q\mapsto q^l$). Hence
	\begin{equation}\label{eq:ch-Schur-functor}
		\ch S_\kappa\big(V[x]\big)=\rho_m\Big(s_\kappa\Big[\frac X{1-q}\Big]\Big),
	\end{equation}
	where $s_\kappa[\frac X{1-q}]=\pi(s_\kappa)$ is the substitution of Appendix~\ref{ssec:plethysm}, since $\pi(p_l)=p_l/(1-q^l)$. Thus the plethysm $X\mapsto X/(1-q)$ is the effect of the functor $V\mapsto V[x]$ on characters. The following theorem identifies the corresponding modules.
	
	\begin{theorem}\label{thm:projectives-plethysm}
		Assume $m>N$ and let $\kappa\vdash N$. The module $S_\kappa(V[x])$ is the projective cover $P_\kappa^{N\Box}$ of $V_\kappa$ in the category $\Rep(\gl_m[x])^{N\Box}$ of graded modules. Consequently
		\[
		\ch P_\kappa^{N\Box}=\rho_m\Big(s_\kappa\Big[\frac X{1-q}\Big]\Big),\qquad P_{\kappa'}^{N_1\Box}\otimes P_{\kappa''}^{N_2\Box}\simeq\bigoplus_{\kappa\vdash N_1+N_2}\big(P^{(N_1+N_2)\Box}_{\kappa}\big)^{\oplus c^\kappa_{\kappa'\kappa''}}.
		\]
	\end{theorem}
	\begin{proof}
		By~\cite{FKhM} (for $m>N$), the module $V[x]^{\otimes N}=\bW_{\omega_1}^{\otimes N}$ is projective in $\Rep(\gl_m[x])^{N\Box}$ and is generated by its zero-degree part $V^{\otimes N}$. The symmetric group acts on $V[x]^{\otimes N}$ by permuting the tensor factors and commutes with $\gl_m[x]$; hence
		\[
		V[x]^{\otimes N}\simeq\bigoplus_{\kappa\vdash N}S_\kappa(V[x])\otimes\mathrm{Sp}_\kappa,
		\]
		with $\mathrm{Sp}_\kappa$ the Specht module. Each $S_\kappa(V[x])$ is a direct summand of a projective module, hence projective; it is generated by its degree-zero component $S_\kappa(V)=V_\kappa$ (a quotient of a module generated in degree $0$ is generated in degree $0$). A graded module generated by a simple degree-zero $\gl_m$-module $V_\kappa$ has the unique maximal graded submodule, the part of positive degree; so $S_\kappa(V[x])$ is projective with the simple top $V_\kappa$, i.e.\ $S_\kappa(V[x])\simeq P_\kappa^{N\Box}$. The character is~\eqref{eq:ch-Schur-functor}, and the tensor product formula follows from the Littlewood--Richardson decomposition of the product of Schur functors, $S_{\kappa'}(U)\otimes S_{\kappa''}(U)\simeq\bigoplus_\kappa S_\kappa(U)^{\oplus c^\kappa_{\kappa'\kappa''}}$, applied to $U=V[x]$.
	\end{proof}
	
	\begin{remark}
		Under Schur--Weyl duality for $\gl_m$ (and its version for the current algebra, \cite[\S3]{KM}) a graded $\mathfrak S_N$-module $M$ corresponds to a graded $\gl_m$-module with the character equal to the Frobenius characteristic of $M$. The module $S^N(V[x])$ corresponds to the polynomial ring $\mathbb{C}[x_1,\dots,x_N]$ with $\mathfrak S_N$ permuting the variables, and its Frobenius characteristic is $h_N[\frac X{1-q}]$ (Lemma~\ref{lem:pleth-h}); similarly $S^N(V[x,y])$ corresponds to $\mathbb{C}[x_1,y_1,\dots,x_N,y_N]$, with the characteristic $h_N\big[\frac X{(1-q)(1-u)}\big]$ ($q$, $u$ record the degrees in $x$, $y$; $u$ is usually denoted $t$ in Haiman's work), and the Garsia--Haiman modules, which are quotients of this ring, are the modules $\GH^{\gl}_\lambda$ of Section~\ref{ssec:GH-fusion}. In Haiman's theory the characters of such modules are expressed by plethysm (the modified Macdonald polynomials $\widetilde H_\lambda(X;q,t)$ are defined through plethystic substitutions); in the present paper the same plethystic substitution is the character of the functor $V\mapsto V[x]$ (or $V\mapsto V[x,y]$): Theorem~\ref{thm:projectives-plethysm} identifies the plethysm of a Schur function with the character of a projective module of the current algebra, and the plethystic substitution inverts it, $\ch M\mapsto\ch M[(1-q)X]$, computing the multiplicities of the projectives in the (infinite) projective resolution of $M$ (Corollary~\ref{cor:Ext-char} below).
	\end{remark}
	
	\subsection{Ext-vanishing and characters}\label{ssec:app-Ext}
	
	We use the structure of highest weight category of~\cite{BBCKL,CI}: the global Weyl modules $\bW_\kappa$ are the standard objects, and the projective $P_\kappa=P^{N\Box}_\kappa$ has a filtration with the subquotients $\bW_\nu$, $\nu\succeq\kappa$, the module $\bW_\kappa$ occurring once, as the top quotient (BGG reciprocity).
	
	\begin{lemma}\label{lem:Ext-Weyl}
		$\Ext^\bullet(\bW_\lambda,V_\kappa)=0$ and $\Ext^\bullet(W_\lambda,V_\kappa)=0$ if $\kappa\not\succeq\lambda$.
	\end{lemma}
	\begin{proof}
		The kernel of the surjection $P_\lambda\to\bW_\lambda$ has a filtration by $\bW_\nu$ with $\nu\succ\lambda$, and applying the same to each of them we obtain a projective resolution of $\bW_\lambda$ consisting of the modules $P_\nu$, $\nu\succeq\lambda$ (the first term is $P_\lambda$). Since $\Hom(P_\nu,V_\kappa)=0$ unless $\nu=\kappa$, we get $\Ext^\bullet(\bW_\lambda,V_\kappa)=0$ if $\kappa\not\succeq\lambda$. By~\S\ref{ssec:app-weyl}, $\bW_\lambda$ is free over the polynomial ring $\bA_\lambda$ and $W_\lambda=\bW_\lambda/\bA_\lambda^+\bW_\lambda$, so the Koszul complex of $\bA_\lambda$ gives a finite resolution of $W_\lambda$ by the modules $\bW_\lambda$ (up to shifts of the grading); the claim for $W_\lambda$ follows.
	\end{proof}
	
	\begin{corollary}\label{cor:Ext-char}
		The following implication holds for finite-dimensional modules $M\in \Rep(\gl_m[x])^{N\Box}$:
		\[
		\left\{\forall \kappa\vdash N,\ \kappa\not\succeq\lambda \quad\Ext^\bullet(M,V_\kappa)=0 \ \right\} \qquad \Rightarrow \qquad
		\ch_qM\in\mathbb Z[q]\Big\langle s_\kappa\Big[\frac X{1-q}\Big]\ :\ \kappa\succeq\lambda\Big\rangle .
		\]
	\end{corollary}
	\begin{proof}
		For $m>N$ the characters $s_\kappa[\frac X{1-q}]$, $\kappa\vdash N$, form a basis of $\Lam^N(q)$ over $\mathbb{C}(q)$, and $\pi$ preserves $\mathbb Z[q]\otimes\Lam$; so $\ch M=\sum_\kappa c_\kappa(q)\,s_\kappa[\frac X{1-q}]$ with $c_\kappa(q)=\langle\ch M[(1-q)X],s_\kappa\rangle\in\mathbb Z[q]$. On the other hand let $P_\bullet\to M$ be the minimal graded projective resolution; its $i$-th term is $\bigoplus_{\kappa,r}P_\kappa\langle r\rangle^{\dim\Ext^i(M,V_\kappa\langle r\rangle)}$, and the generators of the $i$-th term live in degrees $\ge i$, so $\ch M=\sum_i(-1)^i\ch P_i$ in $\mathbb Z[[q]]\otimes\Lam$, i.e.\ $c_\kappa(q)=\sum_{i,r}(-1)^iq^r\dim\Ext^i(M,V_\kappa\langle r\rangle)$. The expansion in the basis is unique, so $c_\kappa=0$ for $\kappa\not\succeq\lambda$.
	\end{proof}
	

\end{document}